\documentclass[preprint]{elsarticle}
\pdfoutput=1
\usepackage[breaklinks=true]{hyperref}
\usepackage{multicol}
\usepackage{amsmath, amsthm, amssymb, dsfont, bm}
\usepackage{cancel}     % Place lines through maths formulae
\usepackage{setspace}
\usepackage{booktabs}
\usepackage{microtype}  % Subliminal refinements towards typographical perfection
\usepackage{graphicx}
\usepackage[table,xcdraw]{xcolor}
\usepackage[american]{babel}
\usepackage{tikz}
\usepackage{pgfplots}
\usepackage[ruled,vlined]{algorithm2e}
\usepackage{soul}

\journal{\ldots}
\graphicspath{{./figures/}}

\newtheorem{remark}{Remark}

\newcommand \kwA {Multi-objective simulation optimization}
\newcommand \kwB {Bayesian optimization}
\newcommand \kwC {Stochastic simulator}
\newcommand \kwD {Gaussian process modeling}
\newcommand \kwE {Kriging}

\hypersetup{
  bookmarksnumbered,
  colorlinks=true,
  pdfauthor={Bruno Barracosa, Julien Bect, Héloïse Dutrieux Baraffe, Juliette Morin, Josselin Fournel, Emmanuel Vazquez},
  pdfkeywords={\kwA, \kwB, \kwC, \kwD, \kwE}
}

\usetikzlibrary{arrows,calc,shapes,positioning,intersections,quotes}
\tikzset{
help lines/.style={dashed, thick},
axis/.style={<->},
important line/.style={thick},
connection/.style={thick, dotted},
}

\usepgfplotslibrary{fillbetween}

\definecolor{mygray}{gray}{0.5}
\definecolor{mysoftgray}{gray}{0.8}

\newcommand{\hlc}[2][yellow]{{\sethlcolor{#1} \hl{#2}}}

\DeclareMathOperator*{\argmax}{argmax}

\newcommand \Rset  {\mathbb{R}}
\newcommand \RR    {\Rset}
\newcommand \Xset  {\mathbb{X}}
\newcommand \XX    {\Xset}
\newcommand \Pc    {\mathcal{P}}
\newcommand \Fc    {\mathcal{F}}
\newcommand \EIm   {\ensuremath{\mathrm{EI}_m}\xspace}
\newcommand \EI    {\ensuremath{\mathrm{EI}} }
\newcommand \PcHat {\widehat{\mskip .5mu\Pc}\mskip -.5mu}
\newcommand \nmax  {n_{\mathrm{max}}}

\begin{document}

\begin{frontmatter}

  \title{Bayesian multi-objective optimization\\ for stochastic
    simulators: an extension of\\ the Pareto Active Learning method}

  \author[EDF,CS]{Bruno Barracosa}

  \author[CS]{Julien Bect}

  \author[EDF]{Héloïse Dutrieux Baraffe}

  \author[EDF]{Juliette Morin}

  \author[EDF]{Josselin Fournel}

  \author[CS]{Emmanuel Vazquez\corref{mycorrespondingauthor}}
  \cortext[mycorrespondingauthor]{Corresponding author}
  \ead{emmanuel.vazquez@centralesupelec.fr}

  \address[EDF]{EDF, Research and Development, Power Systems and Energy Markets, Palaiseau, France}
  \address[CS]{Université Paris-Saclay, CNRS, CentraleSupélec, Laboratoire des signaux et systèmes, 91190, Gif-sur-Yvette, France.}

  \begin{abstract}
    This article focuses on the multi-objective optimization of
    stochastic simulators with high output variance, where the input
    space is finite and the objective functions are expensive to
    evaluate.  We rely on Bayesian optimization algorithms, which use
    probabilistic models to make predictions about the functions to be
    optimized.  The proposed approach is an extension of the Pareto
    Active Learning (PAL) algorithm for the estimation of
    Pareto-optimal solutions that makes it suitable for the stochastic
    setting.  We named it Pareto Active Learning for Stochastic
    Simulators (PALS).  The performance of PALS is assessed through
    numerical experiments over a set of bi-dimensional, bi-objective
    test problems.  PALS exhibits superior performance when compared
    to a selection of alternative approaches based on scalarization
    or random-search.
  \end{abstract}


\end{frontmatter}

\newpage

\section{Introduction}
\label{sec:introduction}

We consider the problem of multi-objective optimization with
stochastic evaluations, also known as multi-objective simulation
optimization \citep{Hunter2019}. %
More precisely, let $f_1, \ldots f_q$ be $q$ real-valued objective
functions defined on a search domain $\mathbb{X} \subset\mathbb{R}^d$,
and assume that these objectives can only be observed with some random
errors: each evaluation of the objectives at a given
point~$x \in \mathbb{X}$ in the search domain yields random responses
$Z_j = f_j(x) + \varepsilon_j$, $1 \le j \le q$, where the
$\varepsilon_j$s are zero-mean random variables (the distribution of
which may depend on~$x$), mutually independent and independent from
past evaluations. %
Our objective is to estimate the Pareto-optimal solutions of the
problem:
\begin{equation}
  \min f_1,\,\ldots,\,f_q\,.
\end{equation}

This setting has many applications in the industry when stochastic
simulators are used to make decisions about uncertain systems with
several, often conflicting, criteria. %
Examples can be found in all areas of engineering and science,
for instance plant breeding \citep{Hunter2016}, aircraft
maintenance \citep{Mattila2014} or electric network planning
\citep{Dutrieux2015}.

In this work, we choose to focus on Bayesian optimization algorithms,
which use probabilistic models to make predictions about the functions
to be optimized. %
One of the classical algorithms for Bayesian \emph{deterministic}
multi-objective optimization is the ParEGO algorithm
\citep{Knowles2006}. %
It relies on a scalarization approach to extend the very popular
single-objective Efficient Global Optimization (EGO) algorithm
\citep{Jones1998}, based on the Expected Improvement (EI) criterion. %
Another scalarization approach, called multi-attribute Knowledge
Gradient \citep{Astudillo2017}, uses the Knowledge Gradient (KG)
criterion instead of the EI criterion. Several Bayesian
multi-objective optimization methods that do not rely on scalarization
have been proposed as well. %
These are, for instance, the well-known Expected Hypervolume
Improvement (EHVI) algorithms \citep{Emmerich2006, Feliot2017}, the
Stepwise Uncertainty Reduction (SUR) approach suggested in
\cite{Picheny2015}, the Predictive Entropy Search for Multi-objective
Bayesian Optimization (PESMO) algorithm \citep{Hernandez2016}, the
Max-value Entropy Search for Multi-objective Bayesian Optimization
(MESMO) \citep{Belakaria2020}, and the Pareto Active Learning (PAL)
algorithm \citep{Zuluaga2013}.

In this article, we choose to focus on PAL because it is easy to
implement and inexpensive in terms of computational time, contrarily
to other Bayesian approaches based on computationally-intensive
criteria for point selection. %
The main idea of PAL consists in using confidence intervals from the
posterior distribution of a Gaussian process to balance exploration
and exploitation, in a way inspired by the Gaussian Process Upper
Confidence Bound rule (GP-UCB) \citep{Srinivas2010} for
single-objective optimization. \cite{Zuluaga2013} present PAL as a
multi-objective Bayesian optimization tool directed at cases where
``evaluating the objective function is expensive and noisy''. %
However, as will be discussed in this article, PAL is fundamentally an algorithm designed for the
optimization of deterministic objective functions over a finite input
space, and significant limitations arise when considering truly
stochastic (random) evaluations.

This article brings two contributions: (1) an extension of the PAL
algorithm, called PALS, making it suitable for stochastic evaluations,
and (2) a numerical benchmark comparing the
performance of the proposed approach with several competing approaches
based on scalarization and random-search.

The article is organized as follows. Section~\ref{sec:mooss}
introduces in more detail the problem of multi-objective optimization
for stochastic simulators. %
Section~\ref{sec:overview} describes the original PAL algorithm. %
Section~\ref{sec:pals} presents the proposed extension for stochastic
simulators. %
Section~\ref{sec:numerical} details the numerical experiments and their
main results. %
Section~\ref{sec:conclu} draws conclusions and suggests several perspectives. %
A separate supplementary material file is available.

\section{Bayesian multi-objective optimization for stochastic simulators}
\label{sec:mooss}

Consider a black-box stochastic simulator whose outputs are noisy
evaluations of $q$ latent functions
$f_1,\ldots,f_q : \Xset \mapsto \Rset$ over a search space
$\Xset \subset \Rset^d$, $d \in \mathbb{N}$. %
The vector of function values $(f_1(x),\ldots,f_q(x))$ will be denoted
by $\boldsymbol{f}(x)$. %
Suppose that the functions $f_1,\,\ldots,f_q$ represent (possibly
conflicting) costs that we want to minimize.

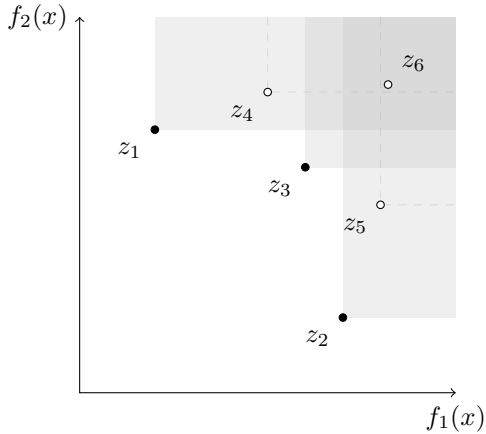
\begin{figure}[tbp]
    \centering
    \begin{tikzpicture}[scale=1]
    %Draw axis
    \coordinate (y) at (0,5);
    \coordinate (x) at (5,0);
    \draw[axis] (y) -- (0,0) --  (x);

    %Coordinates of points
    \coordinate (y1) at (1,3.5);
    \coordinate (y2) at (3.5,1);
    \coordinate (y3) at (3,3);
    \coordinate (y4) at (2.5,4);
    \coordinate (y5) at (4,2.5);
    \coordinate (y6) at (4.1,4.1);
    
            %Domination regions
    \fill[gray!50,nearly transparent] (y1) -- (5,3.5) -- (5,5) -- (1,5) -- cycle;
    \fill[gray!50,nearly transparent] (y2) -- (5,1) -- (5,5) -- (3.5,5) -- cycle;
    \fill[gray!50,nearly transparent] (y3) -- (5,3) -- (5,5) -- (3,5) -- cycle;
    
    % Draw them
    \draw[draw=black,fill=black] (y1) circle [radius=0.05];
    \draw[draw=black,fill=black] (y2) circle [radius=0.05];
    \draw[draw=black,fill=black] (y3) circle [radius=0.05];
    \draw[draw=black,fill=white] (y4) circle [radius=0.05];
    \draw[draw=black,fill=white] (y5) circle [radius=0.05];
    \draw[draw=black,fill=white] (y6) circle [radius=0.05];

    \node[below left=1pt of {y1}, outer sep=1pt,fill=none] {$z_1$};
    \node[below left=1pt of {y2}, outer sep=1pt,fill=none] {$z_2$};
    \node[below left=1pt of {y3}, outer sep=1pt,fill=none] {$z_3$};
    \node[below left=1pt of {y4}, outer sep=1pt,fill=none] {$z_4$};
    \node[below left=1pt of {y5}, outer sep=1pt,fill=none] {$z_5$};
    \node[above right=1pt of {y6}, outer sep=1pt,fill=none] {$z_6$};
    
    \node[below=1pt of {(5,0)}] {$f_1(x)$};
    \node[left=1pt of {(0,5)}] {$f_2(x)$};

    \draw [dashed,gray!80,nearly transparent] (2.5,5) -- (y4) -- (5,4);
    %\draw [dashed,gray!80,nearly transparent] (3,5) -- (y3) -- (5,3);
    \draw [dashed,gray!80,nearly transparent] (4,5) -- (y5) -- (5,2.5);

  \end{tikzpicture}
    \caption{Example of a Pareto front
      $\Fc =\left\lbrace z_1, z_2, z_3 \right\rbrace$, in a
      bi-objective example with finite search
      domain~$\Xset = \{x_1, \ldots, x_6\}$.}
  \label{fig:paretofront}
\end{figure}

In a multi-objective optimization problem, the goal is to identify the
solutions that represent the best possible trade-offs among the
$q$~objectives. %
This is usually defined using the Pareto domination rule: for any two
points $x,\, x^{\prime}\in\RR^{d}$ with images
$z = (z_1,\,\ldots,\, z_q)= \bm{f}(x),\,z^\prime =
(z^{\prime}_1,\,\ldots,\, z^{\prime}_q) = \bm{f}(x^{\prime}) \in
\Rset^q$, we write $z \prec z^\prime$ ($z$ dominates $ z^\prime$ in
the sense of Pareto), when $z_j \leq z^\prime_j$ for all~$j$, with at
least one of the inequalities being strict. %
The set~$\Pc$ of all non-dominated points is called the \emph{Pareto
  set}:
\begin{equation}
  \Pc = \left\lbrace
    x \in \Xset:\; \nexists
    x^\prime \in \Xset,\; \boldsymbol{f}(x^\prime)
    \prec \boldsymbol{f}(x)
  \right\rbrace.
\end{equation}
The \emph{Pareto front} $\Fc \subset \Rset^q$ is by definition the
image $\bm{f}(\Pc)$ of the Pareto set by the objective functions. %
Figure \ref{fig:paretofront} illustrates the Pareto front of a finite
set of evaluations of two objective functions.

The goal of multi-objective optimization is to obtain a good
approximation~$\PcHat_n$ of~$\Pc$ using a sequence
$(X_1,\,\ldots,\, X_n)$, $n\geq 1$, of evaluation points and the
corresponding evaluation results. (The use of capital letters to denote
the evaluation points indicates that these points depend on past
evaluation results.)

In this article, we assume a stochastic setting, where evaluations are
corrupted by additive, normally distributed, homoscedastic
noise. In other words, instead of observing the values of
$f_1,\,\ldots,\, f_q$ at $X_1,\,\ldots,\,X_n$, we observe random
responses
\begin{equation*}
  Z_{i,j} = f_j(X_i) + \varepsilon_{i,j}, \quad 1 \le i \le n,\; 1 \le j \le q,
\end{equation*}
where the random variables $\varepsilon_{i,j}$ are mutually
independent, with Gaussian distributions~$\mathcal{N}(0, \sigma_j^2)$
for some unknown~$\sigma_j^2 > 0$.

Using a Bayesian framework provides a practical solution to handle
this stochastic setting. The main idea is to use a prior probability
distribution, often simply called the prior, to model each objective
$f_j$ as a sample path of a random process~$\xi_j$. %
For convenience, we assume moreover that the $\xi_j$s are independent
Gaussian processes, since it makes it possible to easily obtain the
posterior distributions of the models $\xi_j$ conditional on the
observations $Z_{i,j}$. %
In this case indeed, the $\xi_{j}$s remain independent and Gaussian a
posteriori (i.e., given the evaluation results). %
We denote by $\mu_{n,j}$ and~$k_{n,j}$ their posterior mean and
covariance functions, which can be derived by solving systems of
linear equations \citep[see, e.g.][]{Chiles1999geostatistics,Rasmussen2006}.

\begin{figure}[tbp]
  \centering
  \includegraphics[width=0.7\textwidth]{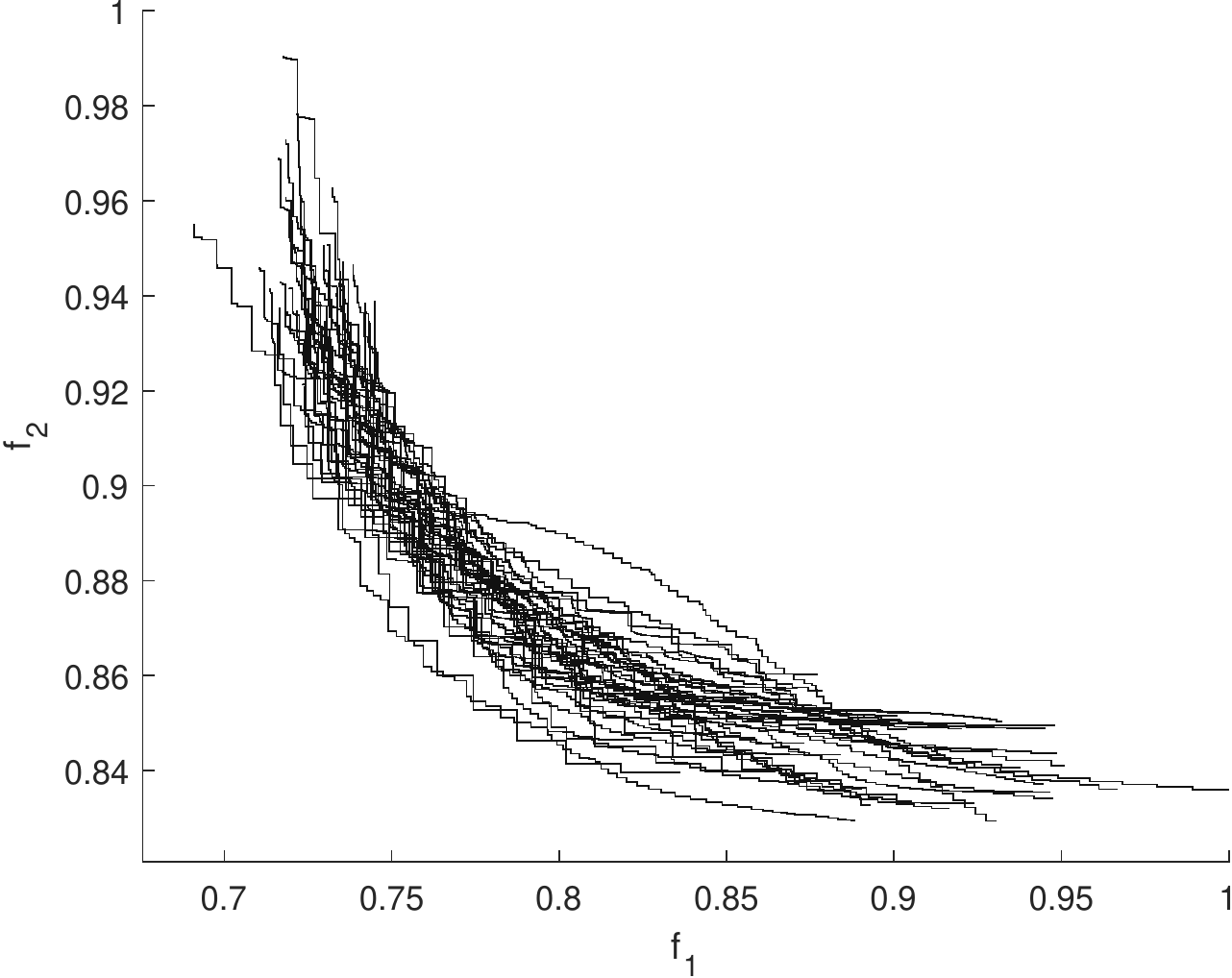}
  \caption{Bi-objective example of normalized random realizations of
    Pareto fronts constructed from probabilistic models.}
  \label{fig:PFIncertitude}
\end{figure}

Posterior distributions based on
past observations can be used simultaneously to choose additional
evaluation points, and to generate predictions of the Pareto front. %
Pareto front predictions can be generated through conditional
simulations (Figure \ref{fig:PFIncertitude}), making it possible not only
to estimate the Pareto front, but also to quantify the remaining
uncertainty \citep{Binois2015}.

As pointed out in the introduction, several Bayesian multi-objective
optimization algorithms have been proposed in the literature, but few
of them can actually be used under a stochastic framework. %
\cite{Astudillo2017} propose an algorithm that takes as input a user
preference and converges to a particular solution of the Pareto
front (i.e., does not provide an estimate of the entire Pareto
set and front). %
The PESMO algorithm \citep{Hernandez2016} on the other hand, which has
been proposed in the machine learning
literature, can be used to jointly estimate the Pareto set and front, as
considered in this article. %
However, this algorithm has a high implementation complexity and
a non-negligible execution cost---it relies at each iteration on
conditional simulations of the GP models, a genetic algorithm and an
expectation-propagation algorithm.

Unlike PESMO, the PAL algorithm proposed by \citet{Zuluaga2013} is a
simple algorithm, but it is not as such suitable for the stochastic
case. %
In the next sections, we introduce the original PAL algorithm, discuss
its limitations, and propose modifications to deal with the stochastic
case.

\begin{remark}
  It is worth noting that several algorithms from the literature of
  budget allocation / ranking and selection, namely MOCBA
  \citep{Lee2010}, SCORE \citep{Pasupathy2014}, and M-MOBA
  \citep{Branke2015} algorithms, use similar ideas to those outlined
  above, but with one notable difference: they assume no correlation
  between the values of the objective functions at distinct points in
  the search space. %
  Recently, an algorithm using a mix of Bayesian optimization and
  ranking and selection approaches was proposed by
  \citet{RojasGonzalez2020} and named SK-MOCBA.
\end{remark}

%%%%%%%%%%%%%%%%%%%%%%%%%%%%%%%%%%%%%%%%%%%%%%%%%%%%%%%%%%%%%%%%%%%%%%%%%%%
%%%%%%%%%%%%%%%%%%%%%%%%%%%%%%%%%%%%%%%%%%%%%%%%%%%%%%%%%%%%%%%%%%%%%%%%%%%
\section{The original PAL algorithm}
\label{sec:overview}

In the following, we assume that $\XX \subset \RR^{d}$ is a (possibly
large) finite set, corresponding for instance to a discrete grid or a
large random sample from a uniform distribution on a bounded set. %
The PAL algorithm builds a sequence~$(X_n)_{n \ge 1}$ of evaluation
points using, at each iteration, a classification of all the points
$x \in \Xset$ of the search space as potentially Pareto-optimal or
not, based on the posterior distributions of the Gaussian
processes~$\xi_j$. %
More precisely, the algorithm partitions~$\Xset$ into three sets at
iteration~$n$: %
a set~$P_n$ of points that are deemed Pareto-optimal, a set~$N_n$
of points that are deemed dominated,
and the set $U_n = \Xset \setminus (P_n \cup N_n)$ of points that remain
unclassified. %
Then, the PAL algorithm improves the approximation of the Pareto set
by choosing an additional evaluation point in~$P_n \cup U_{n}$, where
uncertainty is maximal (in a sense defined below).

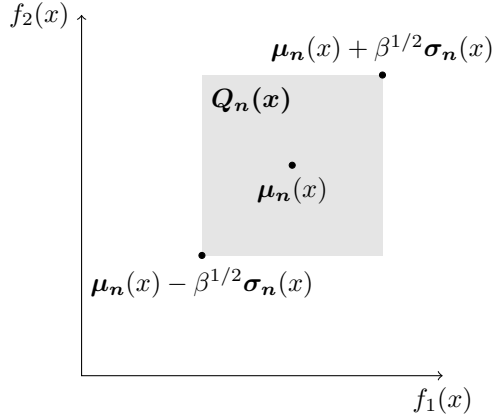
\begin{figure}[tbp]
  \centering
  \begin{tikzpicture}[scale=0.8]
  
    %Draw axis
    \coordinate (y) at (0,6);
    \coordinate (x) at (6,0);
    \draw[axis] (y) -- (0,0) --  (x); 

    \node[below=1pt of {(6,0)}] {$f_1(x)$};
    \node[left=1pt of {(0,6)}] {$f_2(x)$};
    
    %Coordinates of points
    \coordinate (y1) at (2,2);
    \coordinate (y2) at (2,5);
    \coordinate (y3) at (5,5);
    \coordinate (y4) at (5,2);
    
    \coordinate (yc) at (3.5,3.5); %center

	% Box
    %\draw (y1) -- (y2) -- (y3) -- (y4) -- cycle;
    \fill[gray!80,nearly transparent] (y1) -- (y2) -- (y3) -- (y4) -- cycle;
    
    % Draw them
    \draw[draw=black,fill=black] (y1) circle [radius=0.05];
    \draw[draw=black,fill=black] (y3) circle [radius=0.05];
    \draw[draw=black,fill=black] (yc) circle [radius=0.05];

    \node[below=1pt of {y1}] {$\boldsymbol{\mu_n}(x)-{\beta^{1/2} \boldsymbol{\sigma_n}(x)}$};
    \node[above=1pt of {y3}] {$\boldsymbol{\mu_n}(x)+{\beta^{1/2} \boldsymbol{\sigma_n}(x)}$};    
    \node[below=1pt of {yc}] {$\boldsymbol{\mu_n}(x)$};
    \node[below right=1pt and 0pt of {y2}] {$\boldsymbol{Q_n(x)}$};

%    %Domination regions
%    \fill[cyan!50,nearly transparent] (y1) -- (5,3.5) -- (5,5) -- (1,5) -- cycle;
%    \fill[cyan!50,nearly transparent] (y2) -- (5,1) -- (5,5) -- (3.5,5) -- cycle;
%    \fill[cyan!50,nearly transparent] (y3) -- (5,3) -- (5,5) -- (3,5) -- cycle;
%    
%    % Draw them
%    \draw[draw=red,fill=red] (y1) circle [radius=0.025];
%    \draw[draw=red,fill=red] (y2) circle [radius=0.025];
%    \draw[draw=red,fill=red] (y3) circle [radius=0.025];
%    \draw[fill] (y4) circle [radius=0.025];
%    \draw[fill] (y5) circle [radius=0.025];
%    \draw[fill] (y6) circle [radius=0.025];
%
%    \node[below left=1pt of {y1}, outer sep=1pt,fill=none] {$y_1$};
%    \node[below left=1pt of {y2}, outer sep=1pt,fill=none] {$y_2$};
%    \node[below left=1pt of {y3}, outer sep=1pt,fill=none] {$y_3$};
%    \node[below left=1pt of {y4}, outer sep=1pt,fill=none] {$y_4$};
%    \node[below left=1pt of {y5}, outer sep=1pt,fill=none] {$y_5$};
%    \node[above right=1pt of {y6}, outer sep=1pt,fill=none] {$y_6$};
%    
%    \node[below=1pt of {(5,0)}] {$f_1(x)$};
%    \node[left=1pt of {(0,5)}] {$f_2(x)$};
%    
%
%    
%    \draw [dashed,gray!80,nearly transparent] (2.5,5) -- (y4) -- (5,4);
%    %\draw [dashed,gray!80,nearly transparent] (3,5) -- (y3) -- (5,3);
%    \draw [dashed,gray!80,nearly transparent] (4,5) -- (y5) -- (5,2.5);

  \end{tikzpicture}
  \caption{Bi-objective example of a prediction uncertainty region $Q_{n}(x)$.}
  \label{fig:pal_q}
\end{figure}

To build $P_n$, $U_n$ and $N_n$, each $x$ of the search space is
associated to a region $R_n(x)\subset \RR^{q}$, $n\geq 0$, in the objective space,
which quantifies the uncertainty about the values of the objectives at
that point. %
For each $n \geq 0$, the region $R_n(x)$ is constructed as follows. %
First, for each~$x\in\XX$ and each~$j \in \{1, \ldots, q\}$, compute
the posterior means $\mu_{n,j}(x)$ and posterior variances
$\sigma_{n,j}^{2}(x) = k_{n,j}(x,x)$ of the Gaussian random variables
$\xi_j(x)$ modeling the unknown values of the $f_j$s at $x$, given
observations $Z_{1,j},\ldots, Z_{n,j}$ %
(for $n=0$, they are equal to the prior means and variances). %
Then, define hyper-rectangles $Q_{n}(x)\subset\RR^q$ capturing the
uncertainty of the predictions, as illustrated in
Figure~\ref{fig:pal_q}:
\begin{equation}
  Q_{n}(x) = \left\lbrace z \in \Rset^q:
    \boldsymbol{\mu_n}(x)-{\beta^{1/2}\boldsymbol{\sigma_n}(x)}
    \prec z \prec \boldsymbol{\mu_n}(x) + {\beta^{1/2}\boldsymbol{\sigma_n}(x)}
  \right\rbrace\,,
\end{equation}
where
$\boldsymbol{\mu_{n}}(x) = \left( \mu_{n,1}(x), \ldots, \mu_{n,q}(x) \right)$,
$\boldsymbol{\sigma_n}(x) = \left( \sigma_{n,1}(x), \ldots,
\sigma_{n,q}(x) \right)$, and $\beta > 0$ is a positive scaling factor. %
Finally, for each $x\in\XX$, the regions used for classification are
defined as
\begin{equation}
  \label{eq:Rn}
  R_{n}(x) = \left\{
    \begin{array}{ll}
      \{\mu_{n}(x)\} & \text{if } x \in \{X_1, \ldots, X_n\}\,,\\
      R_{n-1}(x) \cap Q_{n}(x) & \text{otherwise},
    \end{array}
\right.
\end{equation}
with the convention that~$R_{-1}(x) = \RR^q$. %
(NB: in the original article by \cite{Zuluaga2013}, a point $x\in\XX$
already visited is artificially assigned a zero posterior variance;
here, we prefer to state that the corresponding region $R_n(x)$
collapses to $\mu_n(x)$, which is equivalent.)

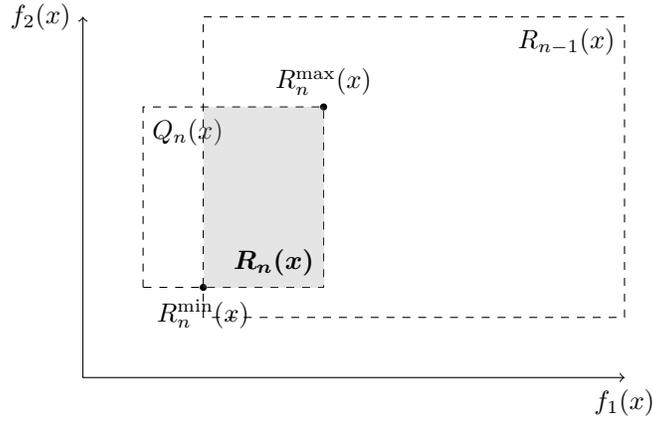
\begin{figure}[tbp]
  \centering
    \begin{tikzpicture}[scale=0.8]
  
    %Draw axis
    \coordinate (y) at (0,6);
    \coordinate (x) at (9,0);
    \draw[axis] (y) -- (0,0) --  (x); 

    \node[below=1pt of {x}] {$f_1(x)$};
    \node[left=1pt of {y}] {$f_2(x)$};
    
    %Coordinates of points
    \coordinate (y1) at (1,1.5);
    \coordinate (y2) at (1,4.5);
    \coordinate (y3) at (4,4.5);
    \coordinate (y4) at (4,1.5);
    
    \coordinate (yc) at (3,3); %center

	% Box
    \draw[dashed] (y1) -- (y2) -- (y3) -- (y4) -- cycle;
    
    \coordinate (y15) at (2,1.5);
    \coordinate (y26) at (2,4.5);
        
    % Draw them
    \draw[draw=black,fill=black] (y15) circle [radius=0.05];
    \draw[draw=black,fill=black] (y3) circle [radius=0.05];

    \node[below=1pt of {y15}] {$R_n^{\min}(x)$};
    \node[above=1pt of {y3}] {$R_n^{\max}(x)$};    
    \node[below right=1pt and 0pt of {y2}] {$Q_n(x)$};
	
    %Coordinates of points
    \coordinate (y5) at (2,1);
    \coordinate (y6) at (2,6);
    \coordinate (y7) at (9,6);
    \coordinate (y8) at (9,1);

		% Box
    \draw[dashed] (y5) -- (y6) -- (y7) -- (y8) -- cycle;
    
    	\fill[gray!80,nearly transparent] (y15) -- (y26) -- (y3) -- (y4) -- cycle;
    	
        \node[below left=1pt and 0pt of {y7}] {$R_{n-1}(x)$};
                \node[above left=1pt and 0pt of {y4}] {$\boldsymbol{R_{n}(x)}$};

  \end{tikzpicture}
  \caption{Bi-objective example of the construction of a
    classification region $R_{n}(x)$ (filled region), and the
    respective optimistic $R_n^{\min}(x)$ and pessimistic
    $R_n^{\max}(x)$ outcomes.}
  \label{fig:pal_intersection}
\end{figure}

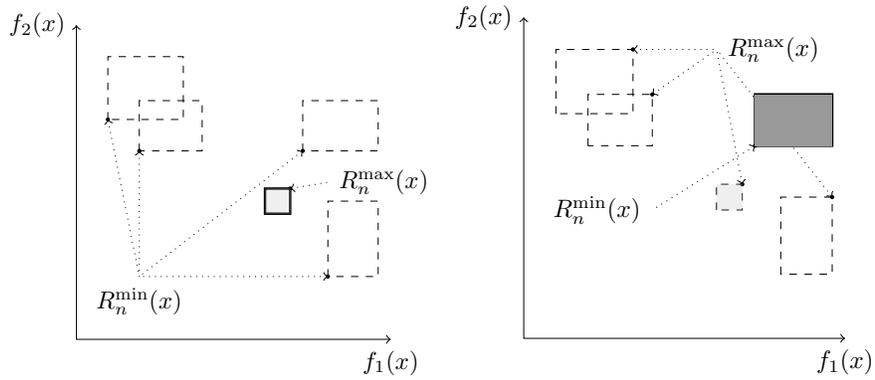
\begin{figure}[tbp]
  \begin{minipage}[t]{0.45\textwidth}
    \resizebox{59mm}{!}{ \begin{tikzpicture}[scale=0.9]
  
    %Draw axis
    \coordinate (y) at (0,5);
    \coordinate (x) at (5,0);
    \draw[axis] (y) -- (0,0) --  (x); 

    \node[below=1pt of {(5,0)}] {$f_1(x)$};
    \node[left=1pt of {(0,5)}] {$f_2(x)$};
    
	% Boxes
	\draw[dashed] (0.5,3.5) -- (0.5,4.5) -- (1.7,4.5) -- (1.7,3.5) -- cycle;
	\draw[dashed] (1,3) -- (1,3.8) -- (2,3.8) -- (2,3) -- cycle;
	\draw [line width=1pt] (3,2) -- (3,2.4) -- (3.4,2.4) -- (3.4,2) -- cycle;
	\draw[dashed] (4,1) -- (4,2.2) -- (4.8,2.2) -- (4.8,1) -- cycle;
	\draw[dashed] (3.6,3) -- (3.6,3.8) -- (4.8,3.8) -- (4.8,3) -- cycle;

    %Coordinates of points
    \coordinate (y1) at (0.5,3.5);
    \coordinate (y2) at (1,3);
    \coordinate (y3) at (4,1);
    \coordinate (y4) at (3.6,3);
    
    \coordinate (y5) at (3.4,2.4);
    
    \coordinate (t1) at (1,1);
    \coordinate (t2) at (4,2.5);
    
    % Arrows
	\draw[<-,dotted] (y1) -- (t1);
	\draw[<-,dotted] (y2) -- (t1);  
	\draw[<-,dotted] (y3) -- (t1);  
	\draw[<-,dotted] (y4) -- (t1);  
	
	\draw[<-,dotted] (y5) -- (t2); 
	
	% Text
	\node[below=1pt of {t1}] {$R_n^{\min}(x)$};
	\node[right=1pt of {t2}] {$R_n^{\max}(x)$};
    
    % Draw them
    \draw[draw=black,fill=black] (y1) circle [radius=0.025];
    \draw[draw=black,fill=black] (y2) circle [radius=0.025];
    \draw[draw=black,fill=black] (y3) circle [radius=0.025];
    \draw[draw=black,fill=black] (y4) circle [radius=0.025];
    
    \draw[draw=black,fill=black] (y5) circle [radius=0.025];
    
    \fill[gray!50,nearly transparent] (3,2) -- (3,2.4) -- (3.4,2.4) -- (3.4,2) -- cycle;

  \end{tikzpicture}}
  \end{minipage}
  \quad
  \begin{minipage}[t]{0.45\textwidth}
    \resizebox{59mm}{!}{  \begin{tikzpicture}[scale=0.9]
  
    %Draw axis
    \coordinate (y) at (0,5);
    \coordinate (x) at (5,0);
    \draw[axis] (y) -- (0,0) --  (x); 

    \node[below=1pt of {(5,0)}] {$f_1(x)$};
    \node[left=1pt of {(0,5)}] {$f_2(x)$};
    
	% Boxes
	\draw[dashed] (0.5,3.5) -- (0.5,4.5) -- (1.7,4.5) -- (1.7,3.5) -- cycle;
	\draw[dashed] (1,3) -- (1,3.8) -- (2,3.8) -- (2,3) -- cycle;
	\draw[dashed] (3,2) -- (3,2.4) -- (3.4,2.4) -- (3.4,2) -- cycle;
	\draw[dashed] (4,1) -- (4,2.2) -- (4.8,2.2) -- (4.8,1) -- cycle;
	\draw [line width=1pt] (3.6,3) -- (3.6,3.8) -- (4.8,3.8) -- (4.8,3) -- cycle;

    %Coordinates of points
    \coordinate (y1) at (1.7,4.5);
    \coordinate (y2) at (2,3.8);
    \coordinate (y3) at (4.8,2.2);
    \coordinate (y4) at (3.4,2.4);
    
    \coordinate (y5) at (3.6,3);
    
    \coordinate (t1) at (3,4.5);
    \coordinate (t2) at (2,2);
    
    % Arrows
	\draw[<-,dotted] (y1) -- (t1);
	\draw[<-,dotted] (y2) -- (t1);  
	\draw[<-,dotted] (y3) -- (t1);  
	\draw[<-,dotted] (y4) -- (t1);  
	
	\draw[<-,dotted] (y5) -- (t2); 
	
	% Text
	\node[right=1pt of {t1}] {$R_n^{\max}(x)$};
	\node[left=1pt of {t2}] {$R_n^{\min}(x)$};
    
    % Draw them
    \draw[draw=black,fill=black] (y1) circle [radius=0.025];
    \draw[draw=black,fill=black] (y2) circle [radius=0.025];
    \draw[draw=black,fill=black] (y3) circle [radius=0.025];
    \draw[draw=black,fill=black] (y4) circle [radius=0.025];
    
    \draw[draw=black,fill=black] (y5) circle [radius=0.025];
    
    \fill[gray!50,nearly transparent] (3,2) -- (3,2.4) -- (3.4,2.4) -- (3.4,2) -- cycle;
    \fill[gray!80] (3.6,3) -- (3.6,3.8) -- (4.8,3.8) -- (4.8,3);

  \end{tikzpicture}}
  \end{minipage}
  \caption{Illustration of the classification of a Pareto-optimal
    point (left) presented in light gray, and the classification of a
    non-Pareto-optimal point (right) presented in dark gray, while
    unclassified points remain as non-filled square dashed outlines.}
  \label{fig:dominationrule}
\end{figure}

For each $R_{n}(x)$, define an optimistic outcome $R_n^{\min}(x)$ and
a pessimistic outcome $R_n^{\max}(x)$ as the best possible outcome and
the worst possible outcome within $R_n(x)$ according to the Pareto
domination rule (see Figure~\ref{fig:pal_intersection}). %
Then, each point $x\in\Xset$ is classified in $P_n$, $U_n$ or $N_n$
according to the following rules illustrated in
Figure~\ref{fig:dominationrule}:
\begin{equation}
  \label{eq:rule-P}
  P_n = \left\{x \in \XX
    \bigm|
    \nexists x^\prime \in \Xset,\;
    R_n^{\min}(x^\prime) +\epsilon \prec R_n^{\max}(x) -\epsilon
  \right\},
\end{equation}
where $\epsilon \in \RR^q_{+}$ is a vector of margin parameters
controlling the acceptation in~$P_n$ and~$N_n$. %
In other words, $x \in P_{n}$ when the ($\epsilon$--shifted)
pessimistic outcome~$R_n^{\max}(x)$ is not dominated by the
($\epsilon$--shifted) optimistic outcome~$R_n^{\min}(x^\prime)$ of any
other point~$x^\prime$. %
Conversely,
\begin{equation}
  \label{eq:rule-N}
  N_n = \left\{x \in \XX
    \bigm|
    \exists x^\prime \in \Xset\setminus\left\lbrace x \right\rbrace,\;
    R_n^{\max}(x^\prime) - \epsilon \prec R_n^{\min}(x) +\epsilon
  \right\},
\end{equation}
or in other words, $x \in N_{n}$ when the ($\epsilon$--shifted)
optimistic outcome~$R_n^{\min}(x)$ is dominated by at least one of the
the ($\epsilon$--shifted) pessimistic outcome~$R_n^{\max}(x^\prime)$
of another point~$x^\prime$; and finally
$U_n = \XX \setminus (P_n \cup N_n)$.

According to \cite{Zuluaga2013}, the margin parameter $\epsilon$
provides flexibility around the uncertainty region, and can be defined
globally or individually for each objective.

The parameter~$\beta$ plays the role of a scaling factor for the
uncertainty regions and is defined at each iteration. %
\cite{Zuluaga2013} suggest using an increasing sequence $\beta_{n}$
defined at each iteration $n$ as
\begin{equation}
  \beta_n=2 \log(q \vert \mathbb{X} \vert \pi^2 n^2 / 6\delta)\,,
  \label{eq:Beta_original}
\end{equation}
with $\delta \in [0,1]$. %
The proposed approach allows to ensure an upper bound of the
hypervolume error with probability $1-\delta$, by specifying
$\delta$. %
Actual convergence and upper bound on the maximal number of needed
iterations are presented for the case where intersecting regions are
used, and a squared exponential kernel is adopted with parameters that
are not re-estimated at each iteration (refer to Corollary 2 in the
original paper for details).

The last step of the algorithm consists in selecting a new evaluation
point~$X_{n+1}$ among the non-visited points in the union
$P_n \cup U_n$ for which uncertainty is maximal, in the sense of the
Euclidean distance between the pessimistic and the optimistic outcomes
(i.e., the diameter of~$R_n(x)$):
\begin{equation}\label{equ:PAL:select-next-point}
  X_{n+1} = \argmax_{x \in (U_n \cup P_n) \setminus S_n}
  \left\lVert R_n^{\min}(x) - R_n^{\max}(x) \right\rVert_2,
\end{equation}
where $S_n$ denotes the set of visited points.

Algorithm~\ref{algo:PALoriginal} summarizes the PAL algorithm.

\begin{remark} \label{rem:normaliz} %
  It appears from the right-hand side
  of~\eqref{equ:PAL:select-next-point} that PAL expects the objectives
  to be suitably normalized, in order to be comparable. %
  Indeed, the criterion at any given point involves a sum of differences of
  objective values:
  \begin{equation*}
    \left\lVert R_n^{\min}(x) - R_n^{\max}(x) \right\rVert_2^2
    \;=\; \sum_{j=1}^q \left( R_{n,j}^{\min}(x) - R_{n,j}^{\max}(x) \right)^2.
  \end{equation*}
  We will assume for simplicity that such a normalization has been
  made beforehand. %
  When this is not possible, the range of the objective
  functions can be learned on-the-fly as done, e.g.,
  by~\cite{Feliot2017}.
\end{remark}

\begin{remark} \label{rem:doe-init} %
  In practice, in Bayesian optimization, it is advisable to start the
  optimization procedure using a set of $n_{0}$ preliminary
  evaluations, which make it possible to select the parameters of the
  Gaussian process model before sequential selection of points.
  This initial set of results is called ``training set''
  by~\cite{Zuluaga2013}, or more traditionally initial design of
  experiments. %
  In this case, the index~$n$ in Algorithm~\ref{algo:PALoriginal}
  corresponds to the number of simulations performed after the initial
  design, and the posterior means and variances at iteration~$n$ are
  in fact computed with respect to~$n + n_0$ evaluation results.
\end{remark}

\begin{algorithm}[tbp]
  \SetKwInOut{Input}{Input}
  \SetKwInOut{Output}{Output}
  \Input{ input space $\mathbb{X}$; GP priors; $\epsilon_q$; $\beta_n$}
  \Output{ predicted Pareto set $\PcHat$}
  $P_{-1} = \varnothing,\;  N_{-1} = \varnothing,\;  U_{-1} = \mathbb{X},\;  S_{-1} = \varnothing$ \\ %%\tcp*{Classification and Evaluated sets} 
  $R_{-1}(x) = \mathbb{R}^q$ for all $x \in \mathbb{X}$ \\
  $n = 0$, all\_classified = false\\
  \Repeat {$\mathrm{all\_classified}$} {
    Obtain $\boldsymbol{\mu}_n(x)$ and $\boldsymbol{\sigma}_n(x)$ for all $x \in \mathbb{X}$,
    $\left\lbrace \boldsymbol{\mu}_n(X_i) = (Z_{i,1}, \ldots, Z_{i,q}) \texttt{ and } \boldsymbol{\sigma}_n(X_i)=0, \forall X_i \in S_n \right\rbrace$ \\
    $R_n(x) = R_{n-1}(x) \cap Q_{\mu_n,\sigma_n,\beta_n}(x)$ for all $x \in \mathbb{X}$ \\
    $P_n=P_{n-1}$, $N_n=N_{n-1}$, $U_n=\varnothing$ \\
    \ForAll{$x \in U_{n-1}$}{
      \uIf{$\left( \nexists x^\prime \in \mathbb{X}\setminus\left\lbrace x \right\rbrace: R_n^{\min}(x^\prime) +\epsilon \prec R_n^{\max}(x) -\epsilon \right) $}{
        $P_n = P_n \cup \left\lbrace x \right\rbrace $
      }
      \uElseIf{$\left( \exists x^\prime \in \mathbb{X}\setminus\left\lbrace x \right\rbrace: R_i^{\max}(x^\prime) -\epsilon \prec R_i^{\min}(x) +\epsilon \right)$}{
        $N_n = N_n \cup \left\lbrace x \right\rbrace $
      }
      \Else{
        $U_n = U_n \cup \left\lbrace x \right\rbrace $
      }}
    \uIf{$U_n= \varnothing$}{
      all\_classified = true
    }
    \Else{
      Choose $X_{n+1} = \argmax_{x \in (U_n \cup P_n) \setminus S_n} \left\lVert R_n^{\min}(x) - R_n^{\max}(x) \right\rVert_2$ \\
      $S_{n+1} = S_n \cup \left\lbrace X_{n+1} \right\rbrace$\\
      Sample $Z_{n+1,j} = f_j(X_{n+1}) + \varepsilon_{n+1,j}$ for all $j \in \{ 1, \ldots, q \}$\\
      $n = n+1$
    }
  }
  return $\PcHat = P_n$
  \caption{The original PAL algorithm}
  \label{algo:PALoriginal}
\end{algorithm}

%%%%%%%%%%%%%%%%%%%%%%%%%%%%%%%%%%%%%%%%%%%%%%%%%%%%%%%%%%%%%%%%%%%%%%%%%%%
%%%%%%%%%%%%%%%%%%%%%%%%%%%%%%%%%%%%%%%%%%%%%%%%%%%%%%%%%%%%%%%%%%%%%%%%%%%
\section{PALS: modification of the PAL algorithm %
  for stochastic simulators}
\label{sec:pals}

The original PAL algorithm presents limitations that make it
unsuitable for stochastic evaluations. %
To circumvent these limitations, we propose several modifications,
which can in fact be seen as \textit{simplifications} of the original
algorithm. %
The simplifications, together with an extension to batch evaluations,
are detailed below. %
As in PAL, it is assumed that $\Xset$ is a finite set. %
The resulting algorithm, summarized in
Algorithm~\ref{algo:PALmodified}, is called PALS (S for Stochastic).

\paragraph{Modification 1: Taking into account prediction variance at
  visited points} %
By artificially setting the prediction variance at visited points to
zero (equivalently, by collapsing~$R_n(x)$ to~$\mu_n(x)$ as
in~\eqref{eq:Rn}), \cite{Zuluaga2013} prevent the PAL algorithm from
visiting any point of the search space more than once. %
This behavior is also enforced by the explicit requirement
(see~\eqref{equ:PAL:select-next-point}) that $X_{n+1}$ should belong
to the set~$\left( P_n \cup U_n \right) \setminus S_n$, where
$S_n$~denotes the set of visited points. %
In fact, this approach can be seen as using a ``nugget'' component
\citep{Matheron63:_princ} in the GP prior, which would be appropriate
for a very irregular deterministic simulator, but is not for a
stochastic simulator---where replicated simulations at a given point
yield independent and identically distributed random responses. %
To circumvent this limitation, we remove the requirement
that~$X_{n+1}$ should belong to~$\Xset \setminus S_n$ and use the
correct posterior variance at all points, visited or not; %
more explicitly, we replace~\eqref{eq:Rn}
and~\eqref{equ:PAL:select-next-point} with
\begin{equation}
  \label{eq:Rn-v1}
  R_{n}(x) =  R_{n-1}(x) \cap Q_{n}(x) \,.
\end{equation}
and
\begin{equation}\label{equ:PALS:select-next-point}
  X_{n+1} = \argmax_{x \in U_n \cup P_n}
  \left\lVert R_n^{\min}(x) - R_n^{\max}(x) \right\rVert_2.
\end{equation}
As a consequence, PALS can visit a given point several times and thus
the number of iterations of the algorithm is not bounded as in~PAL by
the cardinality of the search domain. %
We introduce a user-defined maximum budget of evaluations~$\nmax$,
which will constitute our main stopping criterion (since it is often
extremely costly, in the problems that we consider, to reach a
situation where all the points are classified in~$P_n$ or in~$N_n$). %
We note, however, that keeping the cardinality of $U_n$ as an
additional stopping criterion is the only reason why the distinction
between $U_n$ and $P_n$ is still relevant.

\paragraph{Modification 2: Removing the intersection of uncertainty regions} %
\cite{Zuluaga2013} define~$R_n$ at each iteration as the intersection
of~$R_{n-1}$ and~$Q_n$ in~\eqref{eq:Rn} in order to get a non-increasing
sequences of regions---a monotonicity property which plays a role
in the proof of their theoretical results. %
However, problematic situations can arise as a result of this
definition, where $Q_n(x)$ is not contained in $R_{n-1}(x)$ for
some~$x$ at some iteration~$n$. %
When this happens, $R_n(x)$ is empty and the subsequent behavior of
the algorithm is not properly defined since the pessimistic and
optimistic outcomes at~$x$, $R_n^{\max}(x)$ and~$R_n^{\min}(x)$, are
not defined. %
In the PALS algorithm we remove intersections to address this issue,
and thus further simplify~\eqref{eq:Rn-v1} to
\begin{equation} \label{eq:Rn-v2}
  R_{n}(x) =  Q_{n}(x).
\end{equation}
As a consequence, the subsets~$P_n$ and~$N_n$ defined
by~\eqref{eq:rule-P}--\eqref{eq:rule-N} loose their monotonicity
property as well, and thus membership of all the points in~$\XX$
to~$P_n$, $N_n$ or~$U_n$ must be re-evaluated at each iteration.

\begin{remark}
  As an alternative, it is possible to keep the original idea of using
  intersections and simply ``correct'' the sets~$R_n(x)$ so that
  $R_n(x)$ always contains~$\mu_n(x)$. %
  This can be achieved by taking~$R_n(x)$ to be the smallest
  hyper-rectangle that contains both~$R_{n-1} \cap Q_n(x)$
  and~$\{ \mu_n(x) \}$. %
  Note that doing so artificially fixes the problem of empty regions,
  but nonetheless destroys the monotonicity property, which was the
  original motivation of~\cite{Zuluaga2013} for introducing
  intersections. %
  For the sake of completeness, this approach will also be assessed in
  our experimental section.
\end{remark}

\paragraph{Modification 3: Batches of simulations} %
When the output variability is large, a single simulation provides a
very noisy evaluation of the quantities of interest, which are the
expected responses~$f_j$ of the simulator, and therefore yields by
itself little progress in the optimization procedure. %
A natural idea to gain more information, when parallel computing
resources are available, is to perform several simulations at the same
point at each iteration of the algorithm. %
An approach of this type has been proposed for instance by
\cite{Binois2019}. %
We denote by~$k$ the size the simulation batches, i.e., the number of
simulations performed at the same point at each iteration.

The number of simulation results available after~$n$ iterations with
batch size~$k$ is equal to~$N = n k$, and thus grows quickly if $k$~is
large (for instance $k = 200$ or more in
Section~\ref{sec:numerical}). %
Such a large number of results poses apparently a challenge for
Gaussian process regression, whose computational complexity
is~$O(N^3)$. %
However, this complexity can be reduced, without any approximation,
since Gaussian process regression with Gaussian noise only uses the
empirical mean and variance of the simulation results at each
point~$x \in \Xset$ where at least one simulation was made \citep[see,
e.g.,][]{Dutrieux2015iago, Binois2018practical}: %
a careful implementation leveraging this fact has a complexity that
scales cubically with the number of ``visited'' points, which does not
grow with~$k$.

\begin{remark}
  Another approach, not pursued in this article, would consist in
  sampling simultaneously multiple locations at each iteration, as
  proposed by \cite{Habib2016}, for instance.
\end{remark}

\begin{figure}[tbp]
  \centering
  \includegraphics[width=0.7\textwidth]{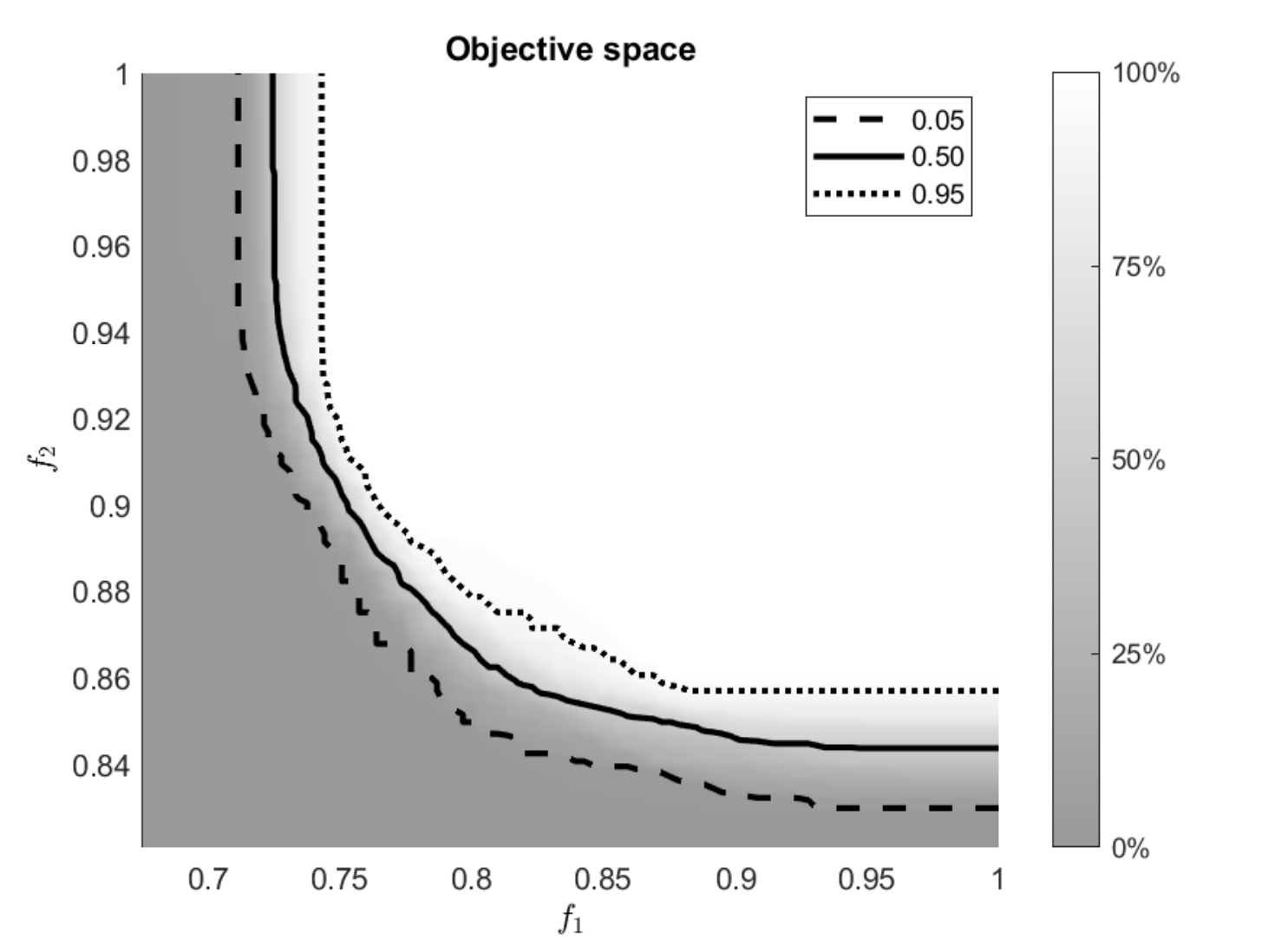}\\[5mm]
  \includegraphics[width=0.63\textwidth]{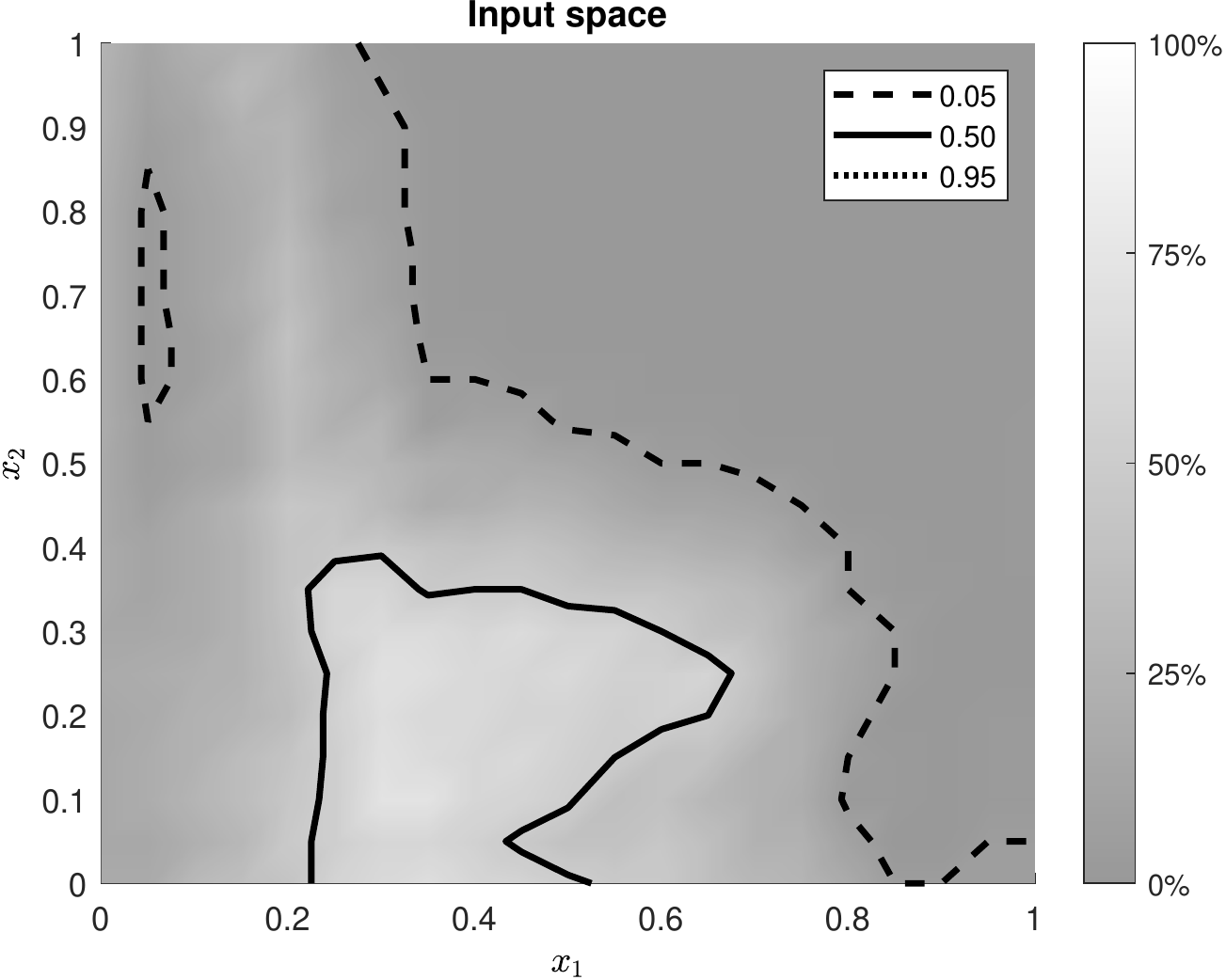}
  \caption{%
    Bi-objective and bi-dimensional example of a normalized heat map
    representing the probability of a point being dominated for the
    objective space (top) %
    and the probability of belonging to the Pareto set (bottom) %
    constructed from a probabilistic model (built from the same data
    and model as in Figure~\ref{fig:PFIncertitude}, %
    using 40 sample paths). %
    Level curves shown for $5\%$, $50\%$, and $95\%$ levels.}
  \label{fig:PFIncertitude_heatmap}
\end{figure}

\paragraph{Predicting the Pareto front and set} %
PAL uses the condition $U_n = \varnothing$ as a stopping criterion,
and then estimates the unknown Pareto set using the set~$P_n$ obtained
at the final iteration. %
With the introduction of a budget-based stopping condition, however,
PALS typically stops before all points are classified. %
A first approach to estimate the Pareto front and set consists in
applying Pareto domination rule directly to the posterior means of the
GP models, in a ``plug-in'' manner. %
This approach, however, does not take the residual uncertainty into
account. %
An alternative approach would be to generate conditional sample paths
of the Gaussian processes instead to represent possible realizations
of the objective functions, and then use them to estimate for each
point in the objective space the probabilities of that point being
dominated (also called empirical attainment function
in~\cite{Binois2015}), and for each point in the input space, the
probabilities of belonging to the Pareto set (also called coverage
probabilities). %
An illustration is provided for the Pareto front and set in
Figure~\ref{fig:PFIncertitude_heatmap}. It is worth noting that a high
confidence about the Pareto front location does not imply a high
confidence about the location of the Pareto set, and vice versa.

Although providing enriched information about the Pareto front and
set, this approach entails a higher computational cost. %
Therefore, for the remaining sections, the plug-in approach is used to
compute estimates of the Pareto front and the Pareto set (at the final
iteration, but also at each iteration for the purpose of algorithm
monitoring and, in our benchmarks, for error assessment).

\begin{algorithm}[tbp]
  \SetKwInOut{Input}{Input}
  \SetKwInOut{Output}{Output}
  \Input{ input space $\mathbb{X}$; GP priors $\mu_{0,j},\sigma_{0,j}, k_j, \forall j$; $\epsilon_j$; $\beta_i$;{\hlc[lightgray]{$n_{\max}$};\hlc[lightgray]{$k$}}.}
  \Output{ predicted Pareto set $\PcHat$.}
  $n = 0$, all\_classified = false\\
  \Repeat {$\mathrm{all\_classified}$ \textbf{or} {\hlc[lightgray]{$n \geq \nmax$}}} {
    {\hlc[lightgray]{Obtain $\boldsymbol{\mu}_n(x)$ and $\boldsymbol{\sigma}_n(x)$ for all $x \in \mathbb{X}$}}\\

    {\hlc[lightgray]{$R_n(x) = Q_{\mu_n,\sigma_n,\beta_n}(x)$ for all $x \in \mathbb{X}$}}

    {\hlc[lightgray]{$P_n=\varnothing$, $N_n=\varnothing$, $U_n=\varnothing$}} \\
    \ForAll{{\hlc[lightgray]{$x \in \mathbb{X}$}}}{
      \uIf{$\left( \nexists x^\prime \in \mathbb{X}\setminus\left\lbrace x \right\rbrace: R_n^{\min}(x^\prime) +\epsilon \prec R_n^{\max}(x) -\epsilon \right) $}{
        $P_n = P_n \cup \left\lbrace x \right\rbrace $
      }
      \uElseIf{$\left( \exists x^\prime \in \mathbb{X}\setminus\left\lbrace x \right\rbrace: R_n^{\max}(x^\prime) -\epsilon \prec R_n^{\min}(x) +\epsilon \right)$}{
        $N_n = N_n \cup \left\lbrace x \right\rbrace $
      }
      \Else{
        $U_n = U_n \cup \left\lbrace x \right\rbrace $
      }
    }
    \uIf{$U_n= \varnothing$}{
      all\_classified = true
    }
    \Else{
      Choose\hlc[lightgray]{$X_{n+1} = \argmax_{x \in (U_n \cup P_n)} \left\lVert R_n^{\min}(x) - R_n^{\max}(x) \right\rVert_2$} \\
      Obtain\hlc[lightgray]{$k$ samples}: $Z_{n+1,j}^{(\ell)} = f_j(X_{n+1}) + \varepsilon_{n+1,j}^{(\ell)},\; 1 \le j \le q, \; 1 \le \ell \le k$ \\
      $n = n+1$
    }
  }
  return $\PcHat$
  \caption{The PALS algorithm}
  \label{algo:PALmodified}
\end{algorithm}

%%%%%%%%%%%%%%%%%%%%%%%%%%%%%%%%%%%%%%%%%%%%%%%%%%%%%%%%%%%%%%%%%%%%%%%%%%%
%%%%%%%%%%%%%%%%%%%%%%%%%%%%%%%%%%%%%%%%%%%%%%%%%%%%%%%%%%%%%%%%%%%%%%%%%%%
\section{Numerical experiments}
\label{sec:numerical}

\subsection{Methodology}

In this section, PALS is compared with a set of alternative
optimization methods: %
1) a pure random search approach (PRS), where evaluation points are
drawn from the uniform probability distribution on~$\XX$; %
2) an alternative random search approach, where evaluation points are
drawn from a model-based probability distribution concentrated on
``promising'' regions; and %
3) two scalarization approaches, which are modifications of the ParEGO
algorithm~\citep{Knowles2006}.

More precisely, the alternative random search approach consists in
drawing the next evaluation point from a distribution proportional to
the probability of that point being misclassified, according to the
posterior distributions of the Gaussian processes $\xi_j$.

For the scalarization techniques, we consider two modified versions of
the ParEGO algorithm suitable for a stochastic setting. %
ParEGO relies on the use of \EI as an infill criterion which is
suitable only for deterministic observations. %
The first version, hereafter named ParEGO-\EIm, uses the \EIm criterion
defined in \cite{Vazquez2008}, which considers the Expected Improvement
with respect to the best prediction over the visited points (as seen
used in the implementation in the PESM branch of Spearmint
\citep{Spearmint}). %
The second one, hereafter named ParEGO-KG, uses the Knowledge Gradient
(KG) sampling criterion instead, which is suitable
for a stochastic setting \citep{Frazier2009}.

For the PALS approach, when an increasing $\beta$ is considered
(cf.~Equation~\eqref{eq:Beta_original}), the parameter~$\delta$ is
taken equal to~$0.05$ as in the original article. %
We consider $\epsilon=0$ for all experiments, given that this
parameter has not shown a significant effect in our numerical
experiments (details available in the supplementary material).

Numerical experiments are performed over a set of nine bi-objective
and bi-dimensional test problems subject to additive homoscedastic
Gaussian noise and defined over a finite $21\times 21$ input space
(see~\ref{annex:problems}). %
A total of 200 optimization runs are performed for each problem and
optimization method.

All the methods except~PRS rely on GP modeling: %
the GP priors on the two objectives are independent, with a constant
unknown mean and a Matérn 5/2 covariance function. %
The unknown mean is integrated out using an improper uniform prior
(``ordinary kriging''), and the parameters of the covariance function
are estimated at each iteration by the Restricted Maximum Likelihood
(ReML) method \citep[see, e.g.,][]{Santner:2003:design}.

Parameters are initialized based on an initial design which is built
by choosing among 1000 initializations of 20 randomly selected points
so as to maximize the minimal distance among them. %
The initial observations are obtained by performing 10 evaluations at
each initial point. %
A total budget~$\nmax$ of~50000 evaluations is allowed at each
experiment in addition to the initial $20 \times 10 = 200$ evaluations. %
At each iteration, a batch of $k$~evaluations is performed at the
selected point. %
The Pareto set and front predictions are generated, for all methods,
using the plug-in approach described in Section~\ref{sec:pals}.

The performance of the methods is assessed using two metrics, averaged
over the 200 runs: %
1) the volume of the symmetric difference between the predicted
dominated region and the true dominated region, %
and 2) the misclassification rate, meaning the percentage of points
incorrectly classified in the Pareto set prediction. %
See~\ref{annex:metrics} for details. %
These metrics provide quantitative information on the error between
the prediction and the correct Pareto set or front.

Unless otherwise mentioned, the setup used for experiments is as
follows: %
no intersection considered, a constant value for $\beta$ corresponding
to the 0.50 probability interval, a null value for $\epsilon$, and a
value of 200 for $k$.

\subsection{Results}

Each of the numerical experiments presented in this section has been
carried out on the nine test problems. %
For brevity, only some representative results are presented; the
remaining results are available in the supplementary material.

\paragraph{Influence of the parameter $\beta$ of PALS}

We are interested in testing the usefulness of an increasing $\beta$
value for the algorithm's performance. %
For this purpose, several constant $\beta$ values were chosen to
generate intervals with different coverage probabilities~$p$. %
More precisely, in this context, we define the coverage probability as
the probability that, for any particular objective~$j$, at any
point~$x \in \Xset$, the function value lies within the region
$\mu_j(x) \pm \beta^{1/2} \sigma_j(x)$. %
Since the posterior is Gaussian, this definition implies that
$\sqrt{\beta} = \Phi^{-1}(0.5 + 0.5p)$, where $\Phi$ denotes the
standard normal cumulative distribution function. %
These values are compared with the increasing-$\beta$ approach
proposed in the original PAL algorithm
(cf. Equation~\eqref{eq:Beta_original}).

The choice of $\beta$ impacts the size of the uncertainty region
around the prediction. %
A small $\beta$ allows the algorithm to focus only on the most
promising regions, by quickly excluding those that seem less
interesting. %
On the other hand, a large~$\beta$ favors exploration by taking longer
to exclude points from being visited.

\begin{figure}[tbp]
  \centering
  \includegraphics[width=0.7\textwidth]{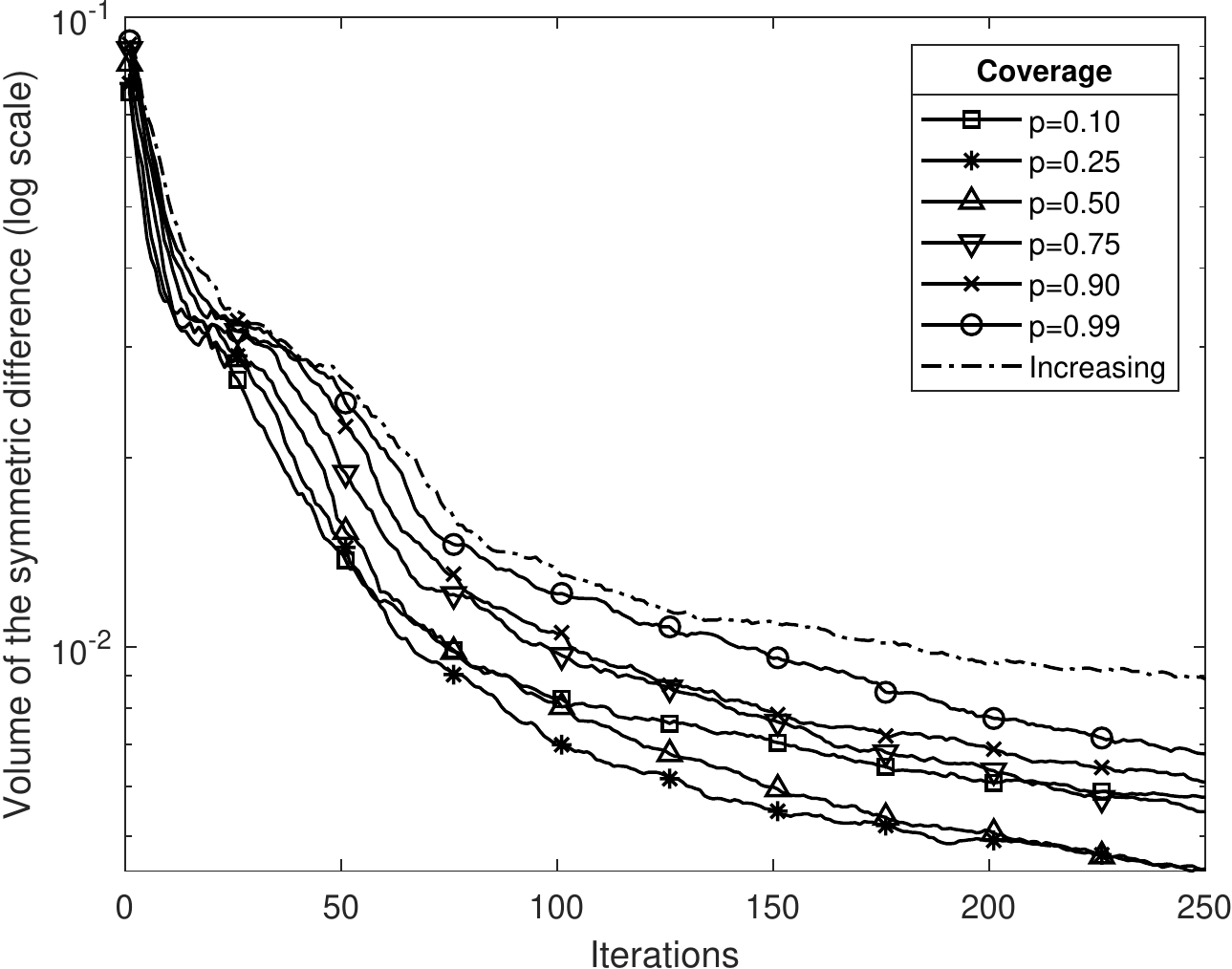}
  \caption{Comparison of average volume of the symmetric difference
    for problem $g_2$, between the increasing-$\beta$ approach and
    fixed $\beta$ for coverage values of 0.10, 0.25, 0.50, 0.75, 0.90,
    and 0.99.}
  \label{fig:BetaHvg4}
\end{figure}

In terms of algorithm performance, an example is presented in Figure
\ref{fig:BetaHvg4} for problem $g_2$, for the volume of the symmetric
difference error metric, comparing the increasing-$\beta$ approach
with fixed $\beta$ for coverage values of 0.10, 0.25, 0.50, 0.75,
0.90, and 0.99.

When looking at the other problem results (as presented in the
supplementary material), the choice of the ideal $\beta$ seems to be
problem-dependent and varies depending on whether the metric used
focuses on the quality of Pareto front or Pareto set predictions. %
However, constant values corresponding to $0.50$ and $0.75$
probability regions have shown an overall interesting performance.

\begin{figure}[tbp]
  \centering
  \includegraphics[width=0.7\textwidth]{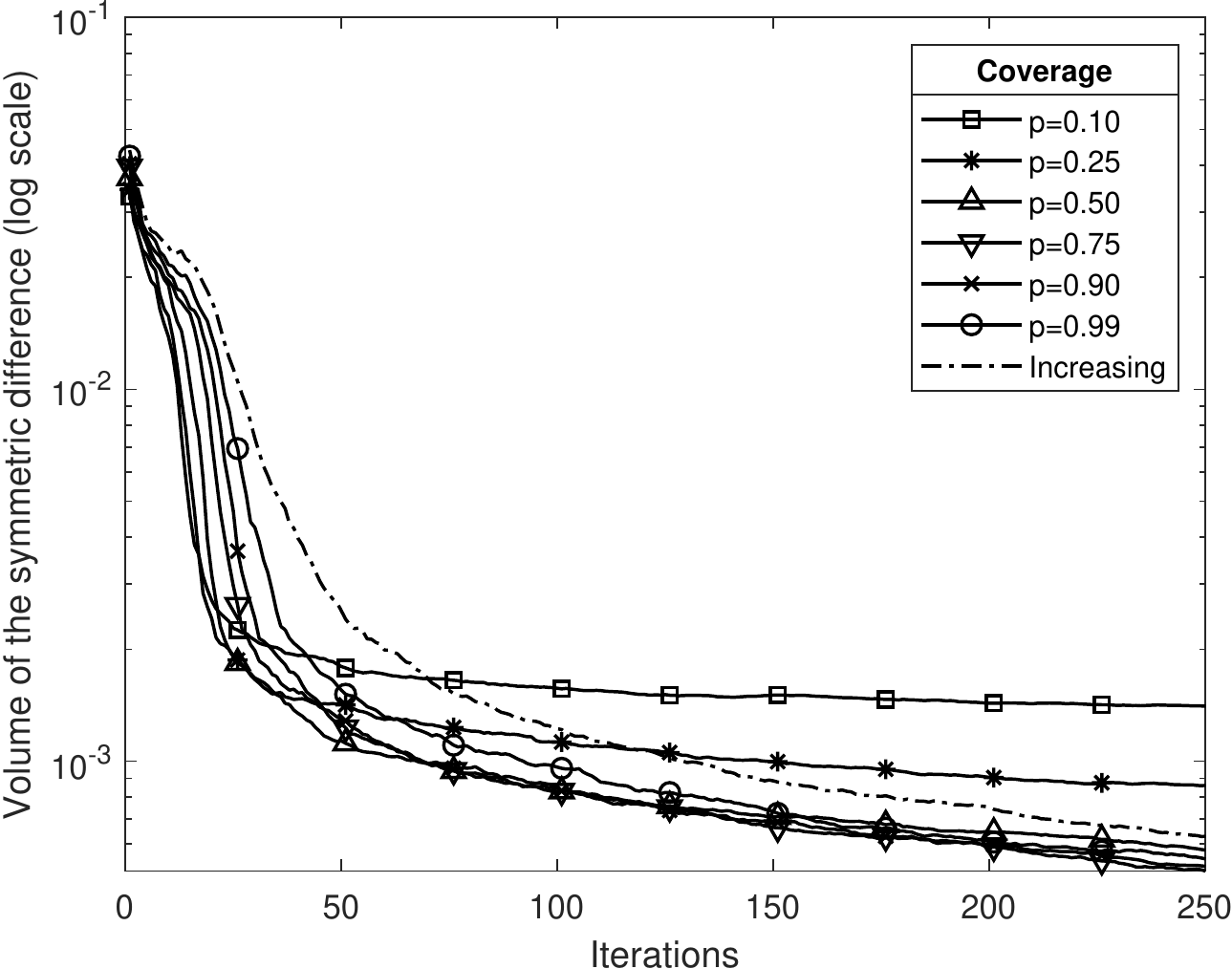}\\[3mm]
  \includegraphics[width=0.7\textwidth]{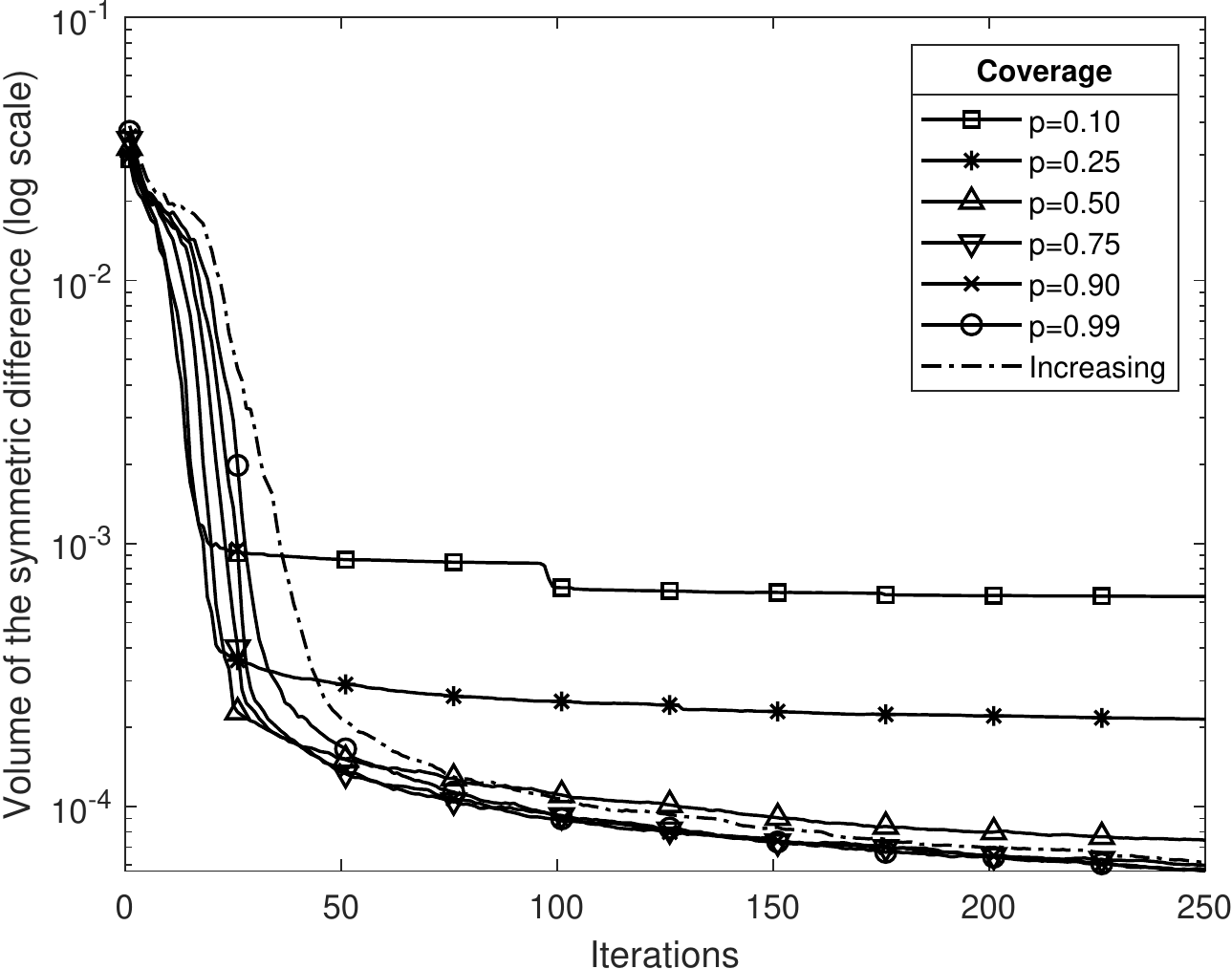}
  \caption{Average volume of the symmetric difference for problem
    $g_2$ for different $\beta$ values, with the problem noise
    standard deviation divided by 10 (top) and by 100 (bottom)}
  \label{fig:BetaHvg4_2}
\end{figure}

Additionally, experiments were performed to assess the impact of the
noise level on the choice of the ideal $\beta$. %
Figure~\ref{fig:BetaHvg4_2} presents the same problem as
Figure~\ref{fig:BetaHvg4} but with two different noise levels. %
In our benchmark, reduced noise levels seem to render the $\beta$
choice less important given the similar performance observed among
choices. %
Small coverage levels (10\% or 25\%), although promising to accelerate
classification, should be avoided given the low performance when
considering the Pareto front prediction quality metric. %
Coverage levels between 50\% and 90\% have shown satisfying results
across a multitude of problems and noise levels for both metrics. %
The increasing-$\beta$ approach (in the form proposed by the PAL
authors) brings no significant advantage when compared to a
fixed~$\beta$ for the classification error metric, but usually
underperforms for the volume of the symmetric difference metric.

Considering the results obtained, a $\beta$ with constant value
corresponding to the $0.50$ probability interval is used in PALS in
the following experiments.

\begin{figure}[tbp]
  \centering
  \begin{minipage}[b]{0.6\textwidth}
    \centering
    \includegraphics[width=\textwidth]{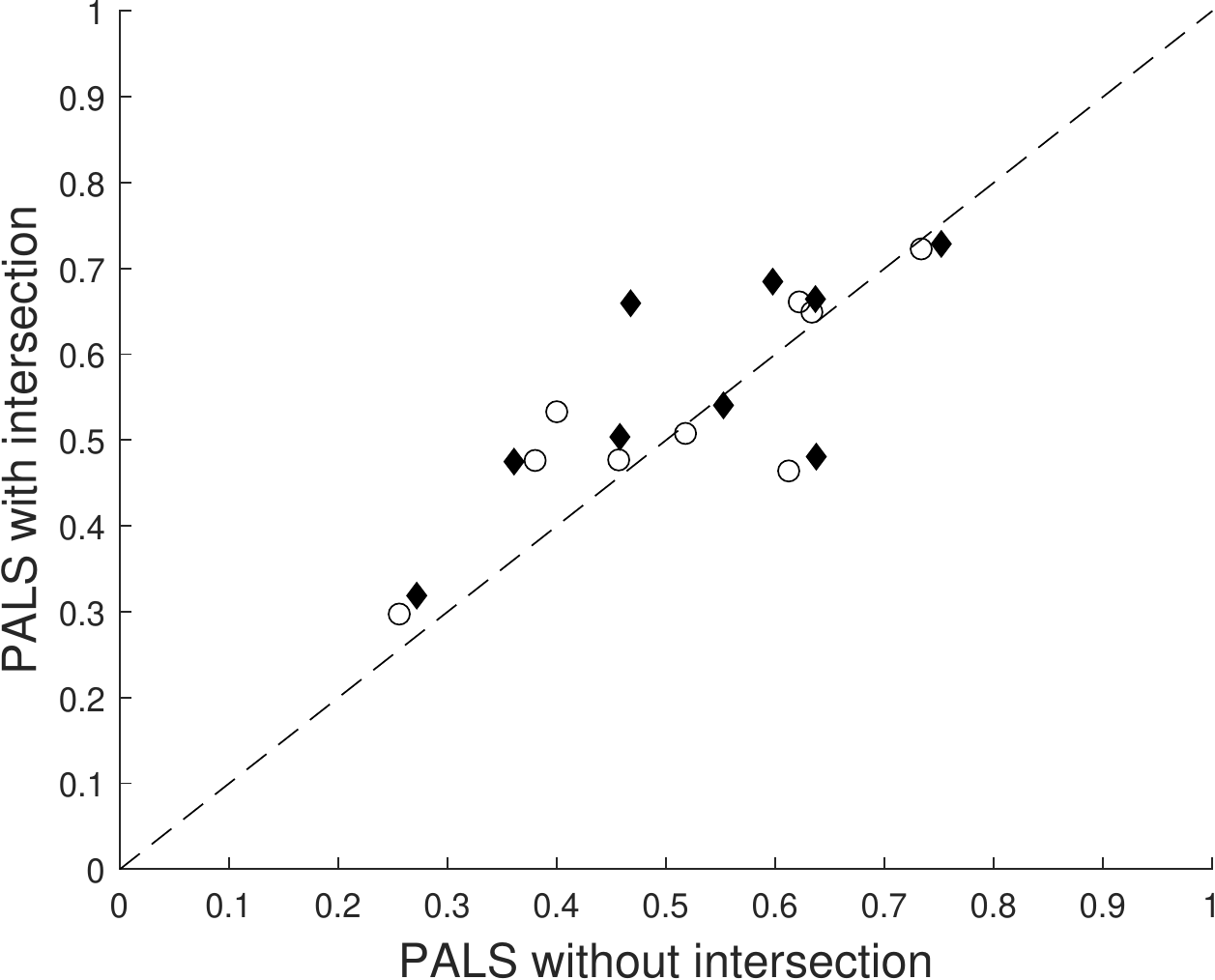}
  \end{minipage}\\[5mm]
  \begin{minipage}[b]{0.6\textwidth}
    \centering
    \includegraphics[width=\textwidth]{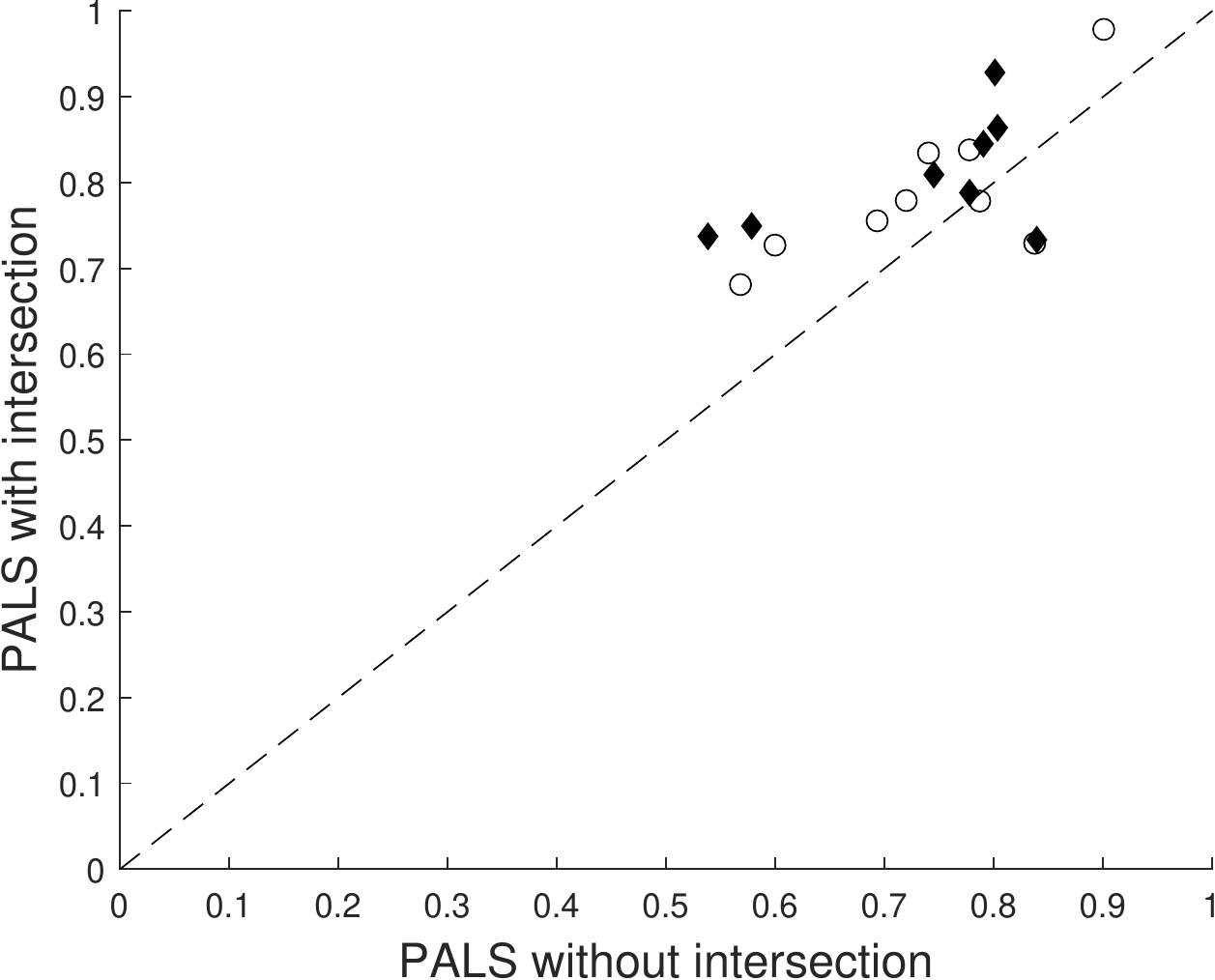}
  \end{minipage}
  \caption{Performance comparison between PALS algorithm with and
    without intersection, normalized by the reference PRS performance.
    Average performance, for each problem for the volume of the
    symmetric difference (top) and misclassification rate (bottom),
    considering 200 random initialization, obtained after a budget of
    25000 evaluations (circles) and 50000 evaluations (diamonds).}
  \label{fig:intersect}
\end{figure}

\paragraph{About intersections} %
In the PALS algorithm, intersecting hyper-regions does not seem to
accelerate or improve the results. %
Figure \ref{fig:intersect} presents a comparison between the
non-intersecting and intersecting variants of PALS by showing for the
different problems the average performance metrics of using (or not)
intersections, normalized by the reference PRS performance. %
It turns out that average results in the middle and final iterations
are very similar, with a slightly better performance of PALS without
intersection for the misclassification rate metric. %
Detailed results can be consulted in supplementary material. %
Given these results, only the non-intersecting version of PALS will be
considered for the remainder of the article.

\begin{figure}[tbp]
  \centering
  \begin{minipage}[b]{0.6\textwidth}
    \centering
    \includegraphics[width=\textwidth]{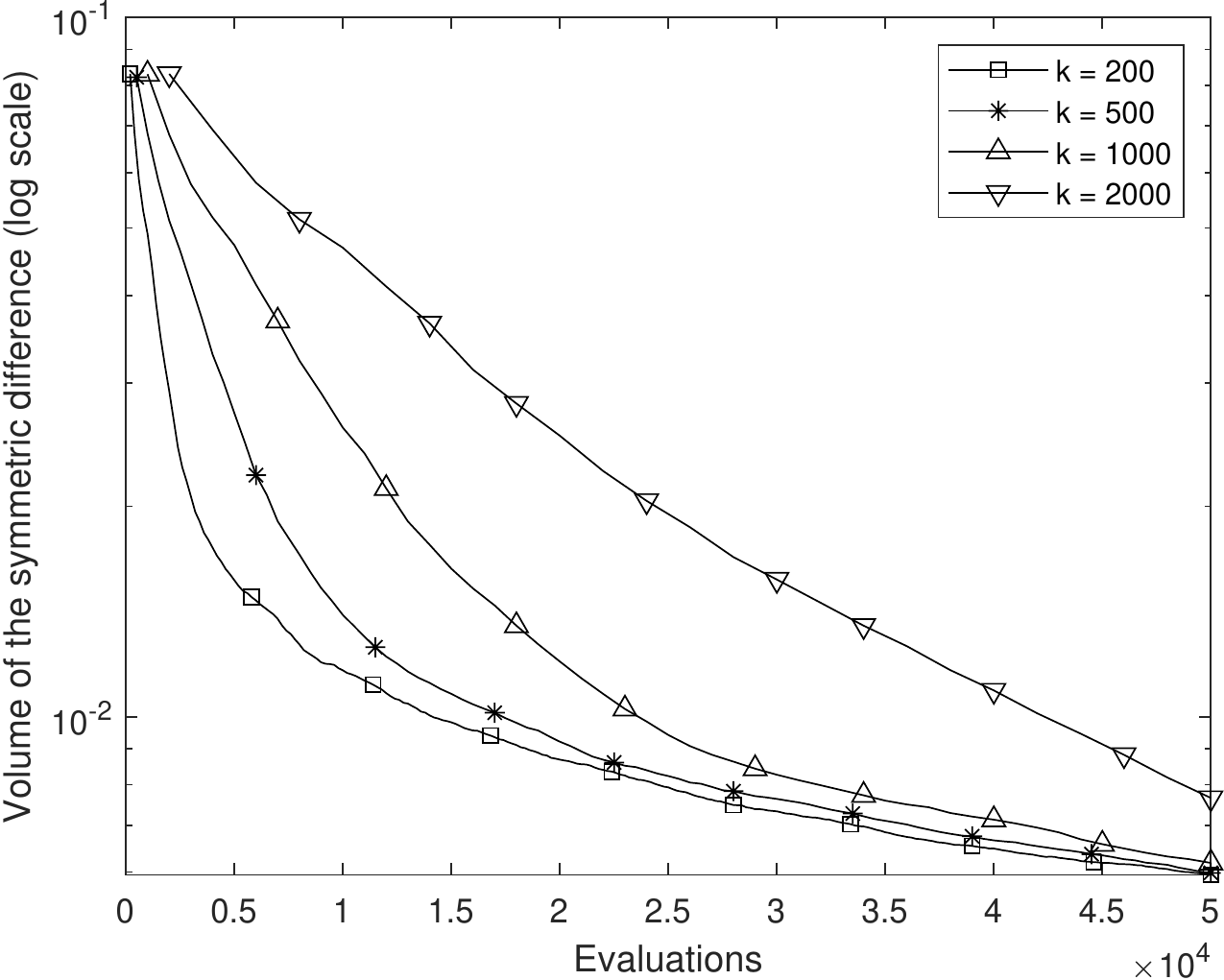}
  \end{minipage}\\[5mm]
  \begin{minipage}[b]{0.6\textwidth}
    \centering
    \includegraphics[width=\textwidth]{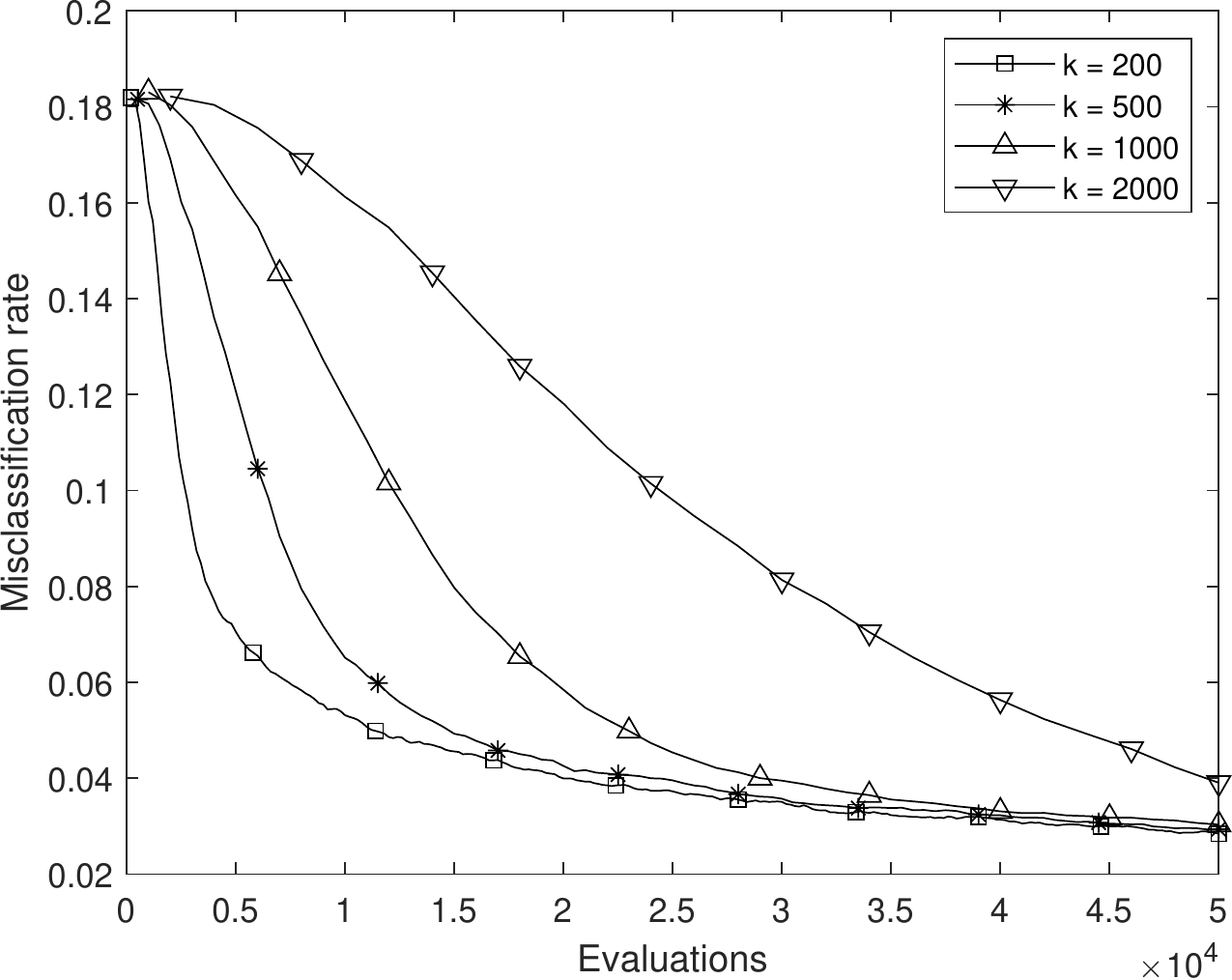}
  \end{minipage}
  \caption{Average performance metrics volume of the symmetric
    difference (top) and misclassification rate (bottom) for problem
    $g_5$, when considering $k$ of 200, 500, 1000 and 2000.}
  \label{fig:kcomparison}
\end{figure}

\paragraph{Influence of the batch size} %
When using a large budget of evaluations, the results show no
significant impact of the batch size $k$ in the average metrics
obtained at the final iteration. %
However, it is worth noting that a smaller batch size for most
problems allowed a faster reduction in metric value at initial
iterations. %
Figure \ref{fig:kcomparison} depicts the comparison between the
average metrics for problem $g_5$ for $k$ values of 200, 500, 1000 and
2000. %
This acceleration of the initial convergence can also be observed for
the other problems on the plots available in the supplementary
material. %
Given the results, $k = 200$ will be used for the remainder of the
article.

\paragraph{Comparison with other approaches} %
This section compares the results of the different approaches: %
PRS, an alternative random search approach (we will name it
``Concentrated Random Sampling'' (CoRS)), two modified ParEGO
approaches, and PALS. %
For PALS, the default parameters are used: a $\beta$ with constant
value corresponding to the 0.50 probability interval, an $\epsilon$ of
0, no intersections, and a $k$ value of 200.

\begin{figure}[tbp]
  \centering
  \begin{minipage}[b]{0.6\textwidth}
    \centering
    \includegraphics[width=\textwidth]{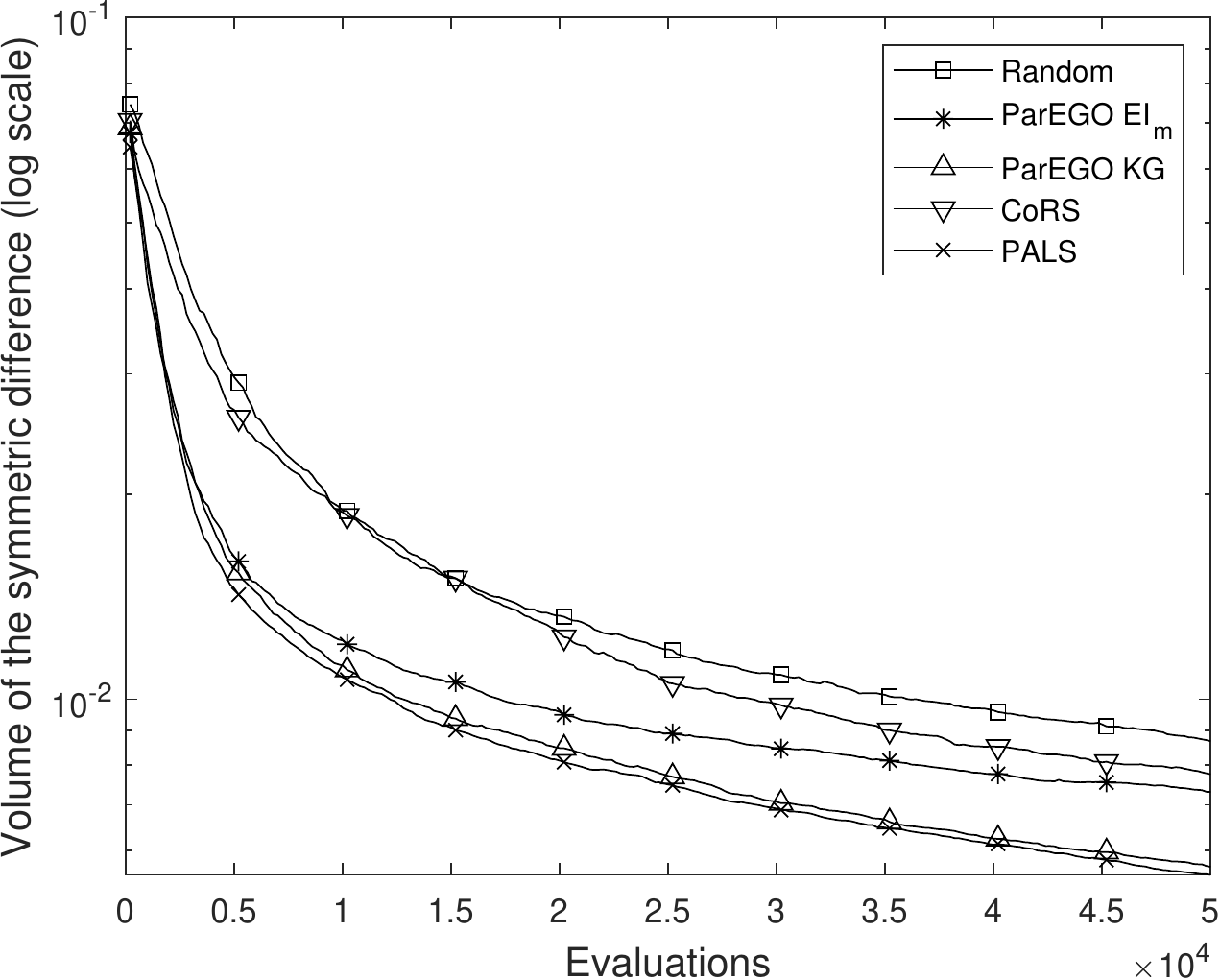}
  \end{minipage}\\[5mm]
  \begin{minipage}[b]{0.6\textwidth}
    \centering
    \includegraphics[width=\textwidth]{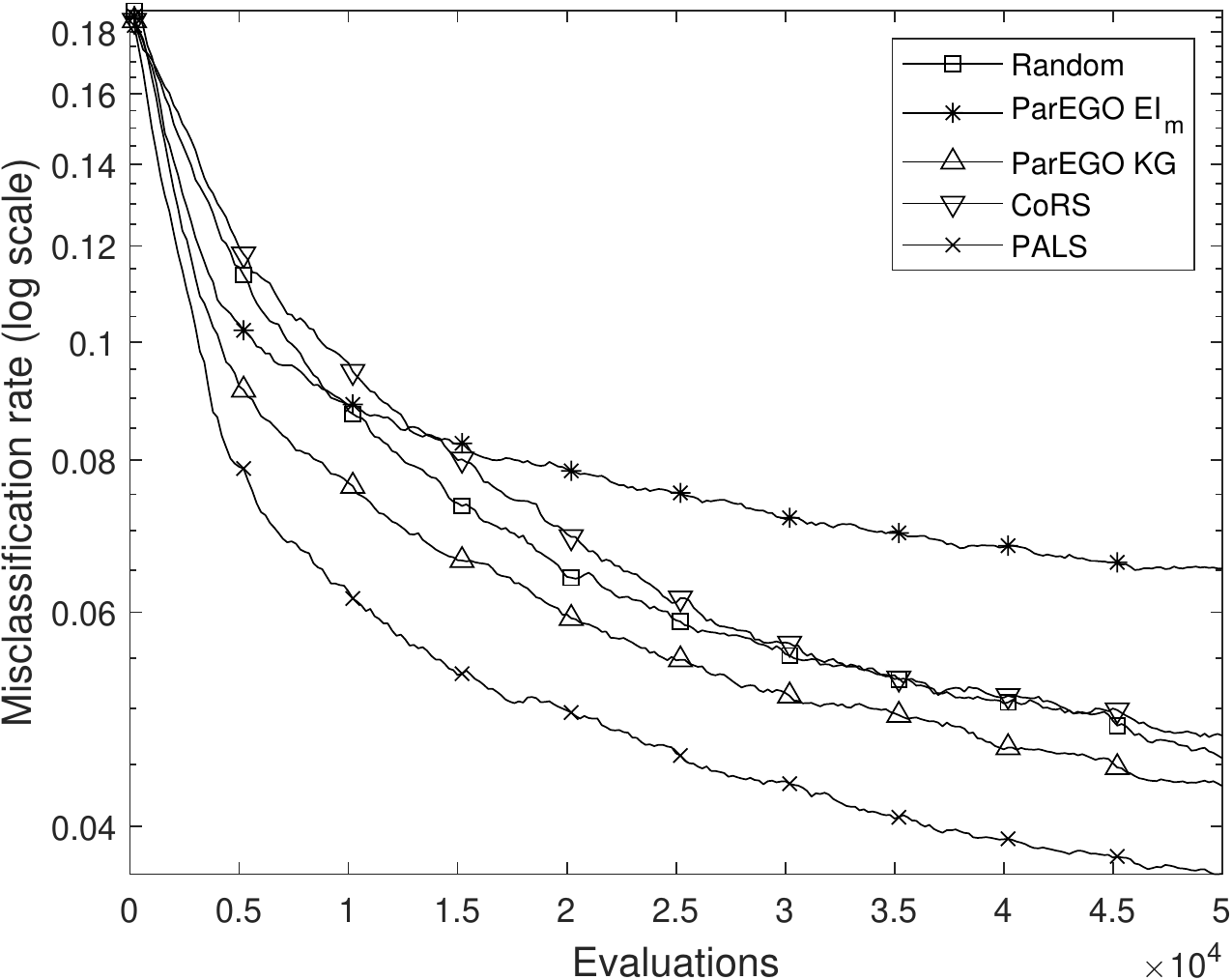}
  \end{minipage}
  \caption{Average volume of the symmetric difference (top) and
    misclassification rate (bottom) on test problem $g_8$ for PRS,
    CoRS, ParEGO-\EIm, ParEGO-KG and PALS.}
  \label{fig:g8metrics}
\end{figure}

Figure \ref{fig:g8metrics} presents the evolution of the performance
metrics on problem $g_8$ for the considered methods. %
Performance illustration for the remaining problems is available in
the supplementary material. %
We can observe that PALS presents a good performance level on both
metrics, allowing for better Pareto set and front estimations in most
cases.

A summary of the numerical results obtained when comparing the
proposed PALS approach and the reference methods over the set of
benchmark problems is presented in Table~\ref{tab:metricresults10}. %
PALS consistently outperforms the other
methods in terms of the ability to better estimate the Pareto set (as
seen through the misclassification rate metric). %
For Pareto front estimation (see volume of the symmetric difference
metric), the most efficient algorithm depends on the problem
considered, with PALS presenting an overall very good performance---although
not being the best performer for problems $g_7$ and $g_9$, it
places second with a metric value close to the leader. %
Additionally, PALS performs significantly better than pure random
search on both metrics for all problems.

\begin{table}[!ht]
  \setlength{\tabcolsep}{3pt}
  \footnotesize
  \centering
  \caption{Average metrics (value in percentage), for the volume of
    the symmetric difference ($V_d$) and misclassification rate ($M$)
    at final iteration, for PRS, CoRS, ParEGO-\EIm, ParEGO-KG, and
    PALS. The best metric values are highlighted in bold with a gray
    background. Metrics at 10\% of the best metric are highlighted
    with a light gray background and considered as acceptable performance. }
  \label{tab:metricresults10}
  \begin{tabular}{@{}crrrrrrrrrrrrrr@{}}
    \cmidrule(l){2-3} \cmidrule(l){5-6} \cmidrule(l){8-9} \cmidrule(l){11-12}\cmidrule(l){14-15}
    \multicolumn{1}{c}{} & \multicolumn{2}{c}{PRS} & & \multicolumn{2}{c}{CoRS} & & \multicolumn{2}{c}{ParEGO-EI$_m$} & & \multicolumn{2}{c}{ParEGO-KG} & &\multicolumn{2}{c}{PALS} \\ 
    \cmidrule(l){2-3} \cmidrule(l){5-6} \cmidrule(l){8-9} \cmidrule(l){11-12} \cmidrule(l){14-15}
    \multicolumn{1}{c}{$g$} & \multicolumn{1}{c}{$V_d$} & \multicolumn{1}{c}{$M$} & & \multicolumn{1}{c}{$V_d$} & \multicolumn{1}{c}{$M$} & & \multicolumn{1}{c}{$V_d$} & \multicolumn{1}{c}{$M$}  & &  \multicolumn{1}{c}{$V_d$} & \multicolumn{1}{c}{$M$} & &  \multicolumn{1}{c}{$V_d$} & \multicolumn{1}{c}{$M$} \\ \midrule
    $g_1$ & 0.793 &7.568 & &0.660 &7.501 & &0.786 &11.418 & &0.687 &\cellcolor{mysoftgray!80}6.211 & &\cellcolor{mygray!80}\textbf{0.596} &\cellcolor{mygray!80}\textbf{6.061} \\ 
    $g_2$ & 0.968 &1.222 & &0.712 &\cellcolor{mysoftgray!80}1.058 & &0.638 &1.330 & &0.502 &1.103 & &\cellcolor{mygray!80}\textbf{0.443} &\cellcolor{mygray!80}\textbf{0.966} \\ 
    $g_3$ & 1.098 &3.491 & &1.143 &3.456 & &\cellcolor{mysoftgray!80}0.707 &3.388 & &\cellcolor{mysoftgray!80}0.711 &\cellcolor{mysoftgray!80}3.019 & &\cellcolor{mygray!80}\textbf{0.700} &\cellcolor{mygray!80}\textbf{2.930} \\ 
    $g_4$ & 1.245 &2.050 & &1.214 &2.180 & &1.025 &2.203 & &\cellcolor{mysoftgray!80}0.756 &1.827 & &\cellcolor{mygray!80}\textbf{0.744} &\cellcolor{mygray!80}\textbf{1.594} \\ 
    $g_5$ & 1.076 &3.815 & &0.749 &3.321 & &0.920 &7.975 & &0.683 &5.743 & &\cellcolor{mygray!80}\textbf{0.594} &\cellcolor{mygray!80}\textbf{2.842} \\ 
    $g_6$ & 1.451 &0.712 & &0.727 &\cellcolor{mysoftgray!80}0.391 & &0.467 &1.545 & &\cellcolor{mysoftgray!80}0.429 &1.049 & &\cellcolor{mygray!80}\textbf{0.394} &\cellcolor{mygray!80}\textbf{0.383} \\ 
    $g_7$ & 0.872 &2.492 & &0.624 &2.539 & &0.512 &4.675 & &\cellcolor{mygray!80}\textbf{0.295} &3.022 & &0.408 &\cellcolor{mygray!80}\textbf{2.230} \\ 
    $g_8$ & 0.868 &4.553 & &0.776 &4.752 & &0.731 &6.517 & &\cellcolor{mysoftgray!80}0.568 &4.320 & &\cellcolor{mygray!80}\textbf{0.552} &\cellcolor{mygray!80}\textbf{3.658} \\ 
    $g_9$ & 1.068 &1.471 & &0.694 &1.169 & &0.639 &2.720 & &\cellcolor{mygray!80}\textbf{0.359} &1.466 & &\cellcolor{mysoftgray!80}0.385 &\cellcolor{mygray!80}\textbf{0.850} \\ 
    \bottomrule
  \end{tabular}
\end{table}

%%%%%%%%%%%%%%%%%%%%%%%%%%%%%%%%%%%%%%%%%%%%%%%%%%%%%%%%%%%%%%%%%%%%%%%%%%%
%%%%%%%%%%%%%%%%%%%%%%%%%%%%%%%%%%%%%%%%%%%%%%%%%%%%%%%%%%%%%%%%%%%%%%%%%%%
\section{Conclusions and future work}
\label{sec:conclu}

In this article, an extension of the Pareto Active Learning (PAL)
algorithm was presented, which makes it suitable for stochastic
simulators with possibly high output variability. %
This new version was named Pareto Active Learning for Stochastic
Simulators (PALS) and showed promising results in a numerical
benchmark on nine bi-objective test problems in dimension two. %
Performance was compared with a Pure Random Search (PRS) approach, a
Concentrated Random Sampling (CoRS) approach and two modified ParEGO
algorithms. %
For the considered benchmark, PALS is the best algorithm to estimate
the Pareto set solutions of the problem. %
It is also the best performer for Pareto front estimation for most of
the problems. %
The modified ParEGO algorithm named ParEGO-KG also
presents promising results when considering the Pareto front
estimation.

We believe that, for future work, comparing the proposed algorithm
with alternative approaches such as SK-MOCBA \citep{RojasGonzalez2020}
or PESMO \citep{Hernandez2016} should be envisaged. %
Additionally, we could be interested in comparisons with algorithms
from the Multi-objective Evolutionary Algorithms (MOEA) literature. %
Interesting perspectives could also include creating hybrid algorithms
that could leverage the different advantages of each algorithm.

A difference in PALS with regards to PAL is that each point selected
during the optimization process is evaluated several times. %
This is especially advantageous when multiple parallel computing units
are available or when the infill criterion to choose the next point to
evaluate is expensive to compute with regards to the actual simulation
time. %
Although the numerical experiments have shown no significant
difference at final iteration between different batch sizes, the
intermediate performances were different. %
Results show that smaller batch sizes have better performance at
initial stages, which seems to indicate a possible interest in having
a noise-level dependent adaptive batch size \citep{Lyu2021adaptive}
during the optimization process.

A main component in the PAL approach is the $\beta$ parameter that
defines the size of the uncertainty region considered around the
model's prediction. %
In the original work, an increasing $\beta$ is proposed. Large
parameter values tend to promote exploration, while smaller values
promote exploitation. %
Numerical experiments have shown, for the problems considered, that
selecting a constant value for the parameter could provide better
results than the original approach. %
Nevertheless, the concept of using an increasing $\beta$ parameter
seems relevant, but further studies should be performed to determine a
suitable initial value and rate of change depending on the problem
structure. %
Therefore, the ideal $\beta$ selection seems to depend on the problem
and noise level, and should still be further studied.

Finally, even though PALS is an algorithm suitable for the case where
the input space is finite, we believe that the main idea could be
applied---with some adaptations---to continuous input spaces.

%%%%%%%%%%%%%%%%%%%%%%%%%%%%%%%%%%%%%%%%%%%%%%%%%%%%%%%%%%%%%%%%%%%%%%%%%%%
%%%%%%%%%%%%%%%%%%%%%%%%%%%%%%%%%%%%%%%%%%%%%%%%%%%%%%%%%%%%%%%%%%%%%%%%%%%
\pdfbookmark{References}{refs}
\bibliographystyle{model5-names} \biboptions{authoryear} % APA style
\bibliography{barracosa-pals-biblio}

\clearpage

%%%%%%%%%%%%%%%%%%%%%%%%%%%%%%%%%%%%%%%%%%%%%%%%%%%%%%%%%%%%%%%%%%%%%%%%%%%
%%%%%%%%%%%%%%%%%%%%%%%%%%%%%%%%%%%%%%%%%%%%%%%%%%%%%%%%%%%%%%%%%%%%%%%%%%%
\pdfbookmark{Appendices}{appendices}
\appendix

%%%%%%%%%%%%%%%%%%%%%%%%%%%%%%%%%%%%%%%%%%%%%%%%%%%%%%%%%%%%%%%%%%%%%%%%%%%
%%%%%%%%%%%%%%%%%%%%%%%%%%%%%%%%%%%%%%%%%%%%%%%%%%%%%%%%%%%%%%%%%%%%%%%%%%%
\section{Benchmark problems}
\setcounter{figure}{0}
\setcounter{table}{0}
\label{annex:problems}

This section introduces the test problems used for the numerical
experiments of Section \ref{sec:numerical}.

All the test functions are presented in Table~\ref{tab:functions} and
are defined in $\left[ 0, 1 \right]^2$. %
Functions~$f_1$ and~$f_2$ mimick a real-life problem presented
by~\cite{Dutrieux2015}. %
Functions $f_3$, $f_4$ (Branin) and $f_5$ (Rosenbrock) are well-known
test functions, rescaled to~$\left[ 0, 1 \right]^2$. %
Functions $f_6$ to $f_{15}$ correspond to randomly generated third
degree polynomials, %
whose coefficients are presented in Table~\ref{tab:coefficients}. %
In our experiments, the values of these 15~test functions are scaled
to~$\left[ 0, 1 \right]$, assuming that the minimum and maximum values
are known (cf.~Remark~\ref{rem:normaliz}).

\begin{table}[tbp]
  \centering
  \caption{Definition of functions $f_1$ to $f_{15}$}
  \label{tab:functions}
  \begin{tabular}{@{}ll@{}}
    \toprule
    $f$ & Definition \\ \midrule
    $f_1$ & $780000 + 110000x_1 - 12000x_2 - 36000x_1x_2 + 280000x_1^2+ 50000x_2^2$ \\
    $f_2$ & $0.83 + 0.17x_1 - 0.015x_2 - 0.0038x_1x_2 + 0.061x_1^2 + 0.0011x_2^2$ \\
    $f_3$ & $e^{(0.36(x_1+x_2))} + 0.6x_1 + 1.2x_2^2 + 3\sin(0.8 \pi x_1)$ \\
    $f_4$ & $(x_2 - a * x_1 .* x_1 + b * x_1 - 6) .^ 2 + c * \cos (x_1) + 10$ \\
    $f_5$ & $100(x_2 - x_1^2)^2 + (1-x_1)^2$ \\
    $f_6$ to $f_{15}$ & $c_1 + c_2x_1 + c_3x_2 + c_4x_1x_2$\\
        & $\quad\mskip -2mu {} + c_5x_1^2 + c_6x_2^2 + c_7x_1^2x_2
          + c_8x_1x_2^2 + c_9x_1^3 + c_{10}x_2^3$ \\ \bottomrule
  \end{tabular}
\end{table}

\begin{table}[tbp]
  \centering
  \caption{Coefficients used for functions $f_6$ to $f_{15}$}
  \label{tab:coefficients}
  \begin{tabular}{@{}lcccccccccc@{}}
    \toprule
    & \multicolumn{10}{c}{Function coefficients} \\ \midrule
    & \multicolumn{1}{l}{$f_6$} & \multicolumn{1}{l}{$f_7$} & \multicolumn{1}{l}{$f_8$} & \multicolumn{1}{l}{$f_9$} & \multicolumn{1}{l}{$f_{10}$} & \multicolumn{1}{l}{$f_{11}$} & \multicolumn{1}{l}{$f_{12}$} & \multicolumn{1}{l}{$f_{13}$} & \multicolumn{1}{l}{$f_{14}$} & \multicolumn{1}{l}{$f_{15}$} \\
$c_1$ & 0.36 & 0.68 & 0.094 & 0.61 & -0.38 & -0.19 & 0.78 & -0.45 & -0.45 & 0.75 \\
    $c_2$ & 8.1 & -9.4 & -7.2 & 5 & 8.5 & 4.8 & 6 & 7.8 & -9.3 & 7.4 \\
    $c_3$ & 7.5 & 9.1 & 7 & 2.3 & 1.4 & 2.1 & -4.7 & -7.7 & -3.5 & -8.2 \\
    $c_4$ & -83 & -2.9 & 49 & -5.3 & 63 & 42 & 90 & 28 & 14 & -98 \\
    $c_5$ & 26 & -60 & 68 & 30 & 81 & 56 & -85 & 34 & -9.7 & 15 \\
    $c_6$ & -80 & 72 & -49 & -66 & 96 & 77 & -82 & -31 & 22 & -31 \\
    $c_7$ & -440 & 160 & 630 & -170 & -120 & 410 & 600 & -500 & -880 & -450 \\
    $c_8$ & 94 & -830 & -510 & -99 & -780 & 360 & 890 & -170 & -370 & -62 \\
    $c_9$ & 920 & -580 & 860 & -830 & -480 & 150 & 370 & -480 & 550 & 780 \\
    $c_{10}$ & 930 & -920 & -300 & 430 & -180 & -16 & -740 & 530 & 390 & -260 \\ \bottomrule
  \end{tabular}
\end{table}

\begin{figure}[p]
  \centering
  \includegraphics[width=0.9\textwidth]{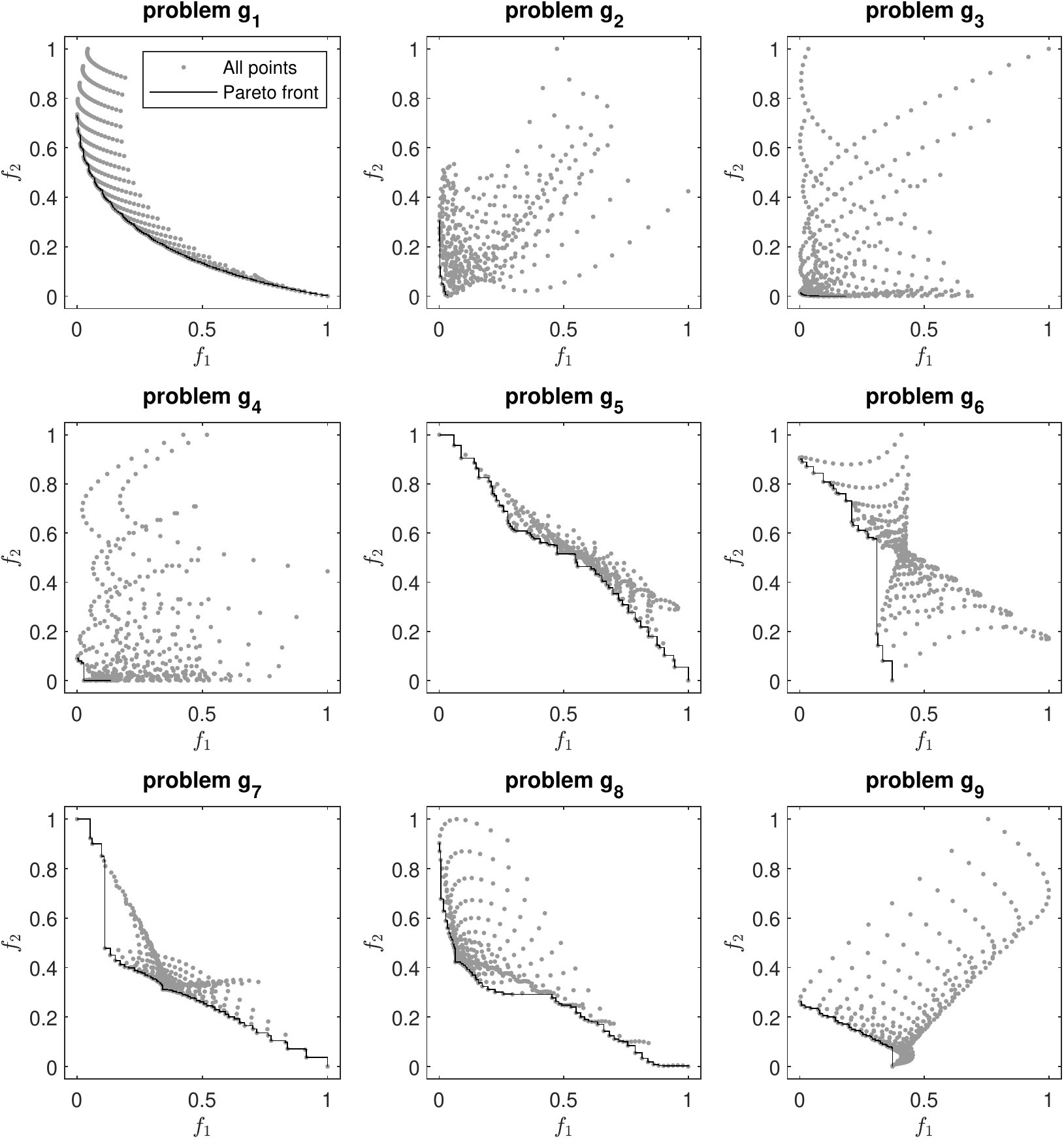}
  \caption{Pareto front solutions for test problems.}
  \label{fig:paretofronts}
\end{figure}

\begin{figure}[p]
  \centering
  \includegraphics[width=0.9\textwidth]{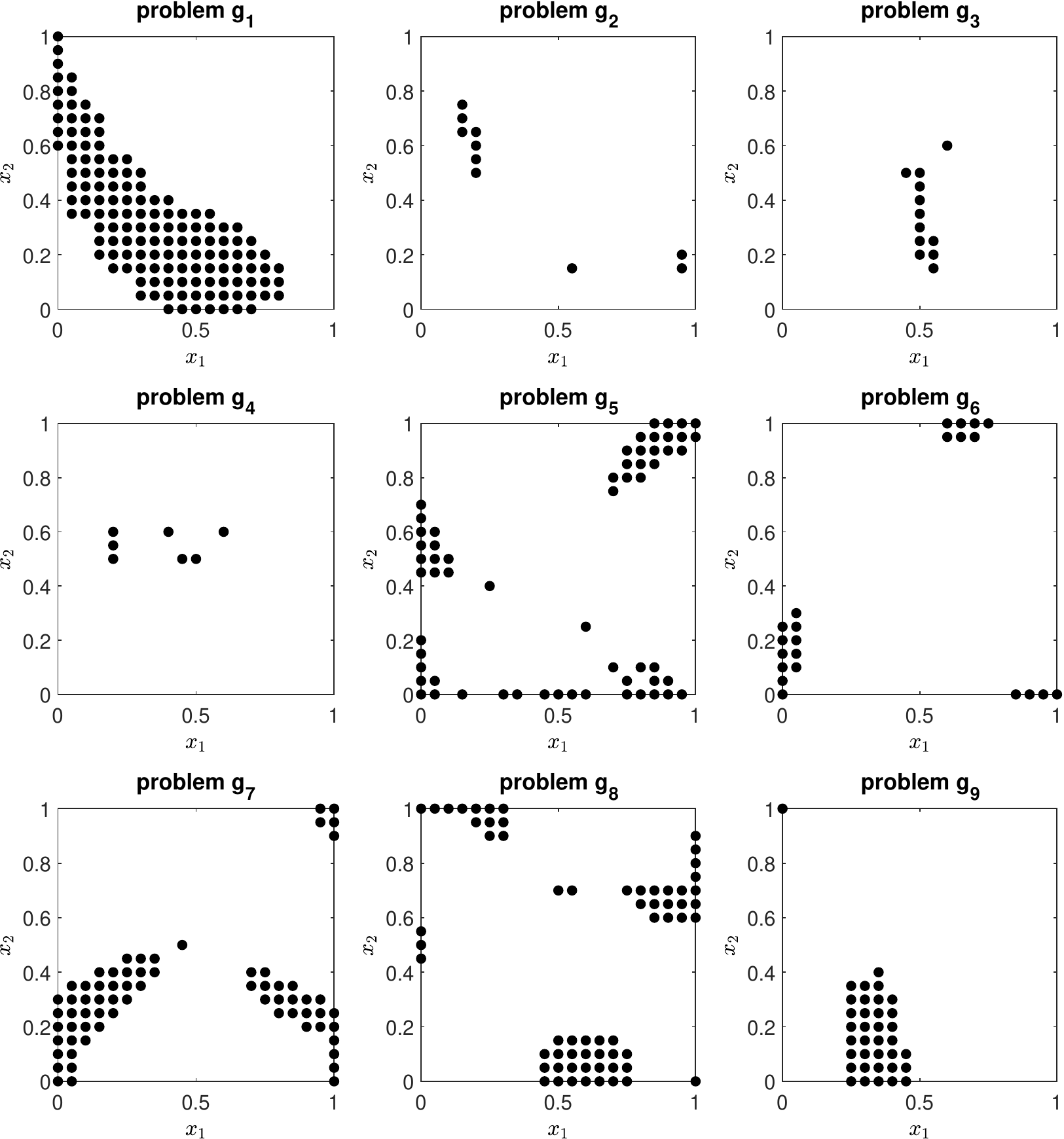}
  \caption{Pareto set solutions for test problems.}
  \label{fig:paretosets}
\end{figure}

A total of nine problems denoted by~$g_1$ to~$g_9$ are proposed based
on the previously defined functions (see Table~\ref{tab:problems}). %
For each problem, the input set~$\Xset$ is a $21 \times 21$ regular
grid on $\left[ 0, 1 \right]^2$. %
The Pareto fronts are presented in Figure~\ref{fig:paretofronts} and
the Pareto sets in Figure~\ref{fig:paretosets}.

\begin{table}[tbp]
  \centering
  \caption{%
    Definition of nine test problems using the test functions
    presented in Tables~\ref{tab:functions}--\ref{tab:coefficients}. %
    For problems~$g_5$ to~$g_9$ the input argument is shifted
    by~$x_0$: if function~$f_k$ is indicated in the table, then the
    actual objective function is $x \mapsto f_k(x - x_0)$. %
    The last two columns provide respectively the cardinality of the
    Pareto set and the proportion of Pareto-optimal solutions in the
    21x21 input set.}
  \label{tab:problems}
  \begin{tabular}{@{}cccllccc@{}}
    \cmidrule(l){2-5}
    \multicolumn{1}{l}{} & \multicolumn{2}{c}{objectives} & \multicolumn{2}{c}{noise variance} \\ \cmidrule(l){2-8} 
$g$ & $f_1^g$ & $f_2^g$ & \multicolumn{1}{c}{$\sigma_1^2$} & \multicolumn{1}{c}{$\sigma_2^2$} & \multicolumn{1}{c}{$x_0$} & $\mid \Pc \mid$ & $\mid \Pc \mid / \mid \mathds{X} \mid$\\ \midrule
    $g_1$ & $f_1$ & $f_2$ & $3.6\times10^9$ & $3.9\times 10^{-3}$ & & 136 & 31\% \\
    $g_2$ & $f_4$ & $f_3$ & $3.1\times10^2$ & $4.8\times10^3$ & & 10 & 2\% \\
    $g_3$ & $f_4$ & $f_5$ & $3.1\times10^2$ & $5.7\times10^8$ & & 12 & 3\% \\
    $g_4$ & $f_3$ & $f_5$ & $4.8\times10^3$ & $5.7\times10^8$ & & 7 & 2\% \\
    $g_5$ & $f_6$ & $f_7$ & $7.0\times10^2$ & $5.6\times10^3$ & $(0.5,0.5)$ & 60 & 14\% \\
    $g_6$ & $f_8$ & $f_9$ & $5.8\times10^2$ & $3.1\times10^3$ & $(0.5,0.5)$ & 22 & 5\% \\
    $g_7$ & $f_{10}$ & $f_{11}$ & $2.1\times10^3$ & $3.2\times10^2$ & $(0.5,0.5)$ & 67 & 15\% \\
    $g_8$ & $f_{12}$ & $f_{13}$ & $1.4\times10^4$ & $1.6\times10^3$ & $(0.3,0.8)(0.6,0.6)$  & 63 & 14\% \\
    $g_9$ & $f_{14}$ & $f_{15}$ & $3.7\times10^3$ & $2.0\times10^4$ & $(0.3,0.8)$ & 36 & 8\% \\ \bottomrule
  \end{tabular}
\end{table}

%%%%%%%%%%%%%%%%%%%%%%%%%%%%%%%%%%%%%%%%%%%%%%%%%%%%%%%%%%%%%%%%%%%%%%%%%%%
%%%%%%%%%%%%%%%%%%%%%%%%%%%%%%%%%%%%%%%%%%%%%%%%%%%%%%%%%%%%%%%%%%%%%%%%%%%
\section{Performance metrics}
\setcounter{figure}{0}
\setcounter{table}{0}
\label{annex:metrics}

We introduce the performance metrics proposed to assess prediction
quality. %
The volume of the symmetric difference measures a global error in the
estimation. %
We use it to assess Pareto front estimates. The misclassification rate
is used as a global error metric to evaluate Pareto set estimates.

\subsection{Volume of the symmetric difference}

\newcommand \RRef {\bm{R}}

Let the ``dominated region'' defined between a Pareto front $\Fc$ and
a reference point~$\RRef$ be the set of all the points that simultaneously
dominate~$\RRef$, and are dominated by a member of $\Fc$:
$D(\Fc) = \cup_{y^\star \in \Fc} \left\lbrace \boldsymbol{y} \in
  \mathbb{R}^q : \boldsymbol{y^\star} \prec \boldsymbol{y} \prec
  \RRef \right\rbrace$ . %
Recalling that our objective functions are all scaled
to~$\left[ 0, 1 \right]$, we take $\RRef = \left( 1.1, 1.1 \right)$
in our (bi-objective) experiments. %
The volume of the dominated region is denoted by $V(\Fc)$, and is illustrated in
Figure \ref{fig:hypervolumedifference} (left).

We can then define the volume of the symmetric difference of the
hypervolumes defined between the true Pareto front and the estimation,
and a reference point~$\RRef$.  We denote it $V_d(\Fc^\star,\hat{\Fc})$:
\begin{equation}
  V_d(\Fc^\star,\hat{\Fc})
  = V(D(\Fc) \setminus D(\Fc^\star)) \cup (D(\Fc^\star) \setminus D(\Fc)),
\end{equation}
and illustrate it in Figure \ref{fig:hypervolumedifference} (right):

\begin{figure}[tbp]
  \centering
  \begin{minipage}[t]{0.45\textwidth}
    \centering
    \resizebox{55mm}{!}{%\documentclass[border=10pt]{standalone}
%\usepackage{pgfplots}

%\usetikzlibrary{arrows,calc,shapes,positioning,intersections,quotes}

%\pgfplotsset{width=7cm, compat=1.10}
%\usepgfplotslibrary{fillbetween}
\pgfmathdeclarefunction{poly}{0}{\pgfmathparse{((x-6)^2)/3}}
\pgfmathdeclarefunction{polyt}{0}{\pgfmathparse{((x-6)^2+2*sin(500*sqrt(x)))/3}}
\pgfmathdeclarefunction{top}{0}{\pgfmathparse{10}}
%\begin{document}
\begin{tikzpicture}[scale=0.8]
  \begin{axis}[
    axis y line = left,
    axis x line = bottom,
    xtick       = {0},
    xticklabels = {},
    ytick       = {0},
    yticklabels = {},
    samples     = 160,
    domain      = 0:5,
    xmin = 0, xmax = 5,
    ymin = 0, ymax = 10,
  ]
  \addplot[name path=poly, black, thick, mark=none, ] {poly};
  \addplot[name path=polyt, black, dashed, mark=none, ] {polyt};
  \addplot[name path=top, white, thick, mark=none, ] {top};
  %\addplot[name path=line, gray, no markers, line width=1pt] {3};
  \addplot fill between[ 
    of = poly and top, 
    split, % calculate segments
    every even segment/.style = {white!70},
    every odd segment/.style  = {gray!60,nearly transparent}
  ];
  
	  \coordinate (R) at (500,100);
  
    \node[below left=1pt of {R}] {$R$};
  
	  %\draw[draw=black,fill=black] (R) circle [radius=3];
  
\end{axis}
	  \draw[draw=black,fill=black] (R) circle [radius=0.05];
\end{tikzpicture}
%\end{document}
}
  \end{minipage}
  \quad
  \begin{minipage}[t]{0.45\textwidth}
    \centering
    \resizebox{55mm}{!}{%\documentclass[border=10pt]{standalone}
%\usepackage{pgfplots}
%\pgfplotsset{width=7cm, compat=1.10}
%\usepgfplotslibrary{fillbetween}
\pgfmathdeclarefunction{poly}{0}{\pgfmathparse{((x-6)^2)/3}}
\pgfmathdeclarefunction{polyt}{0}{\pgfmathparse{((x-6)^2+2*sin(500*sqrt(x)))/3}}
%\begin{document}
\begin{tikzpicture}[scale=0.8]
  \begin{axis}[
    axis y line = left,
    axis x line = bottom,
    xtick       = {0},
    xticklabels = {},
    ytick       = {0},
    yticklabels = {},
    samples     = 160,
    domain      = 0:5,
    xmin = 0, xmax = 5,
    ymin = 0, ymax = 10,
  ]
  \addplot[name path=poly, black, thick, mark=none, ] {poly};
  \addplot[name path=polyt, black, dashed, mark=none, ] {polyt};
  %\addplot[name path=line, gray, no markers, line width=1pt] {3};
  \addplot fill between[ 
    of = poly and polyt, 
    split, % calculate segments
    every even segment/.style = {gray!70,nearly transparent},
    every odd segment/.style  = {gray!70,nearly transparent}
  ];
  
  	  \coordinate (R) at (500,100);
  
    \node[below left=1pt of {R}] {$R$};

\end{axis}
	  \draw[draw=black,fill=black] (R) circle [radius=0.05];
\end{tikzpicture}
%\end{document}
}
  \end{minipage}

  \caption{Hypervolume generated (gray area) between a Pareto front
    prediction $\hat{\Fc}$ (filled line) and a reference point $R$, in
    a bi-objective example (left), and volume of the symmetric
    difference (gray area) between a Pareto front $\Fc$ (dashed line)
    and a prediction $\hat{\Fc}$ (filled line) using reference point
    $R$, in a bi-objective example (right).}
  \label{fig:hypervolumedifference}
\end{figure}
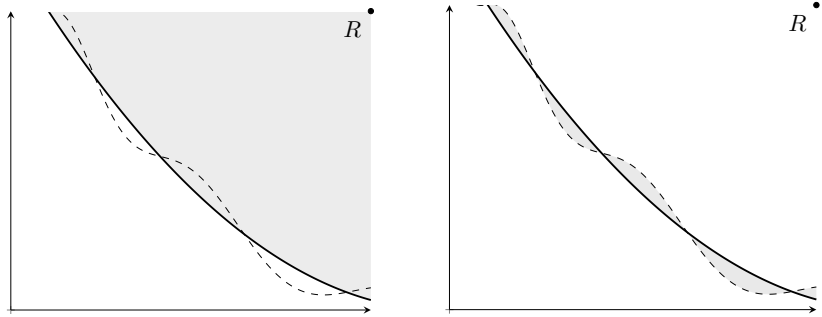

\subsection{Misclassification rate}

The misclassification rate provides a measure of the error of a Pareto
Set prediction. %
Consider a true Pareto set $\Pc^\star$ and a Pareto set prediction
$\PcHat$. %
Then the misclassification rate $M(\Pc^\star,\PcHat)$ is defined as
the proportion of points that differ between the prediction and true
sets:
\begin{equation}
  M(\Pc^\star,\PcHat) =\frac{1}{\lvert \mathbb{X} \rvert}
  \sum_{x \in \mathbb{X}} \left(
    \mathds{1}_{x \in \Pc^\star } - \mathds{1}_{x \in \PcHat}
  \right)^2.
\end{equation}

\end{document}